# Relatively projective groups as absolute Galois groups

JOCHEN KOENIGSMANN

*Dedicated to Yuri Ershov on the occasion of his 60-th birthday*

Submitted to Israel Journal of Mathematics, May 1, 2000

**Abstract**

By two well-known results, one of Ax, one of Lubotzky and van den Dries, a profinite group is projective iff it is isomorphic to the absolute Galois group of a pseudo-algebraically closed field. This paper gives an analogous characterization of relatively projective profinite groups as absolute Galois groups of regularly closed fields.

# Contents

# Introduction

The **absolute Galois group** $G_F$ of a field $F$ is the Galois group of a separable closure $F^{sep}$ of $F$ over $F$, considered as profinite group.

The **free product** $G_1 \star \cdots \star G_n$ of profinite groups $G_1, \ldots, G_n$ is a profinite group $G$ allowing embeddings $\epsilon_i : G_i \to G$ ($i = 1, \ldots, n$) such that, given any homomorphisms $\gamma_i : G_i \to H$ ($i = 1, \ldots, n$) into a profinite group $H$, there is a unique homomorphism $\gamma : G \to H$ with $\gamma_i = \gamma \circ \epsilon_i$ for each $i$.

This paper has two targets: one is to give a simplified proof of the fact that the free product of finitely many absolute Galois groups is again an absolute Galois group, and the other is to describe the absolute Galois group of multiply valued fields satisfying a local-global principle for rational points of varieties.

**Theorem 1** *Given fields $F_1, \ldots, F_n$, there is a field $F$ of characteristic $0$ with $G_F \cong G_{F_1} \star \cdots \star G_{F_n}$. Moreover, if char $F_1 = \ldots =$ char $F_n = p > 0$, $F$ can also be chosen to have characteristic $p$.*

We call a profinite group $G$ **projective [strongly projective] relative to subgroups** $G_1, \ldots, G_n$ **of** $G$ if each epimorphism $\pi : H \twoheadrightarrow G$ of profinite groups which splits locally (i.e. $\forall i \exists \rho_i : G_i \to H$ with $\pi \circ \rho_i = id_{G_i}$) splits globally (i.e. $\exists \rho : G \to H$ with $\pi \circ \rho = id_G$ [and for each $i$, $\rho(G_i)$ is conjugate to $\rho_i(G_i)$ in $H$]).

If $G$ is projective relative to subgroups $G_1, \ldots, G_n$, then $G$ embeds into $G_1 \star \cdots \star G_n \star F$, where $F$ is some free profinite group (Proposition 1.4(5)). Since any free profinite group occurs as absolute Galois group of some field of any prescribed characteristic, and since subgroups of absolute Galois groups are absolute Galois groups, Theorem 1 immediately generalizes to

**Theorem 1'** *Let $G$ be a profinite group which is projective relative to subgroups $G_1, \ldots, G_n$ and assume that each $G_i$ is an absolute Galois group. Then $G$ is an absolute Galois group of some field of characteristic $0$. Moreover, if all $G_i$ can be realized over fields of the same fixed positive characteristic, then so can $G$.*

Theorem 1, which answers Problem 18 from [J], should be attributed to Florian Pop, though he never states it: Theorem 1 is a simple consequence of Theorem 3.4 in [Po1], which even allows to generalize Theorem 1 to certain *infinite* free products of absolute Galois groups. The theorem was first stated and proved quite differently, in [M] in the case where the $G_i$ are of countable rank. Then (without the 'Moreover'), Ershov was the first to publish a proof ([Er3], Theorem 3), which is more in the spirit of Pop's proof.

In all approaches, the technique for realizing free products of given absolute Galois groups is valuation theoretic: find a field $F$ where each of the given absolute Galois groups occurs as decomposition subgroup of $G_F$ w.r.t. some valuation on $F$, and make sure that these valuations are 'in sufficiently general position' to ensure that the decomposition subgroups *freely* generate a subgroup of $G_F$. In [M], this is achieved by 'probabilistic' methods: if these valuations live on a countable Hilbertian field, then, with probability 1, random conjugates of the decomposition subgroups generate a free product ([Ge], Theorem 4.1). This method, however, only works for Galois groups which are isomorphic to subgroups of absolute Galois groups of countable Hilbertian fields, i.e. for countably generated absolute Galois groups. In [Po1] and [Er3], the valuations were put in sufficiently general position by constructing a field which is also *regularly closed* (see below) w.r.t. finitely many valuations having the prescribed decomposition groups.

Many of the arguments in our proof can be found in [Po1] and [Er3], but our proof becomes easier for three reasons: one is that it seems unnecessary to construct a multiply valued field which is regularly closed, the second is that we work with a very handy criterion for profinite groups to be the free product of given subgroups in terms of solving 'locally split embedding problems' (Proposition 1.2). And, thirdly, we restrict ourselves to fields with a finite (rather than a boolean) family of valuations and, thus, avoid all the machinery needed to handle the more general situation.

Theorem 1' generalizes Theorem 3.4 in [Po1] (for finitely many subgroups), since our notion of relative projectivity happens to coincide with Pop's. This coincidence is based on our characterization of relatively projective group (Proposition 1.4) which provides an analogue to Gruenberg's characterization of projective groups ([Gr], Proposition 1).

Our second target is to improve another result of Pop and Ershov on the absolute Galois group of regularly closed fields. Let us recall that an $n$-fold valued field $(F, v_1, \ldots, v_n)$ (with corresponding henselisations $F_1, \ldots, F_n$) is **regularly closed** (or **pseudo-closed**) if it satisfies a **local-global-principle** for rational points on varieties, i.e. if every absolutely irreducible (affine) $F$-variety with a simple $F_i$-rational point for each $i$ has an $F$-rational point.

We shall prove

**Theorem 2** *Let $(F, v_1, \ldots, v_n)$ be a regularly closed n-fold valued field, let $F_1, \ldots, F_n$ be henselisations of $F$ w.r.t. $v_1, \ldots, v_n$ resp., and assume that $v_1, \ldots, v_n$ are independent. Then $G_F$ is strongly projective relative to $G_{F_1}, \ldots, G_{F_n}$.*

Using the variant of Theorem 1' which - as in the proofs of Pop and Ershov - realizes strongly projective groups over *regularly closed* fields (Corollary 5.2), Theorem 2 can immediately be strengthened to the following relative analogue of the Ax-Lubotzky/van den Dries characterization of projective profinite groups as absolute Galois groups of PAC-fields:

**Theorem 2'** *Let $G$ be a profinite group with subgroups $G_1, \ldots, G_n$ where each $G_i$ is isomorphic to some absolute Galois group. Then $G$ is strongly projective relative to $G_1, \ldots, G_n$ iff $G$ is isomorphic to the absolute Galois group of a regularly closed n-fold valued field $(F, v_1, \ldots, v_n)$, where $v_1, \ldots, v_n$ are independent and where, for each $i$, the isomorphism $G \cong G_F$ maps $G_i$ onto a decomposition subgroup of $G_F$ w.r.t. $v_i$.*

Under the hypothesis of our Theorem 2, [Po1], Theorem 3.3 resp. [Po2], Theorem 3.2, comes to the weaker conclusion that $G_F$ is projective relative to $G_{F_1}, \ldots, G_{F_n}$, and only under additional hypotheses (e.g. that the $G_i$ be isomorphic to absolute Galois groups of real or $p$-adically closed fields) strong projectivity has been proved ([Po1], Theorem 1.2). Erhov shows Theorem 2 ([Er2], Thm. 3) for what he calls '$RC^\star$-fields', which are special regularly closed fields, but not all regularly closed fields are $RC^\star$-fields. Both Pop's and Ershov's results, however, deal with certain infinite families of valuations, not just finite ones. We will extend our Theorem 2 to this setting later.

What made our progress on Pop's and Ershov's achievements possible, is the observation that solvability of finite embedding problems for absolute Galois groups is an *existential* first-order property in the language of fields (Observation 5.3), and that the same holds for the corresponding 'relative embedding problems' in the language of $n$-fold valued fields (Proposition 5.4). This is based on a careful analysis of decomposition subfields of finite Galois extensions which goes beyond usual ramification theory (Lemma 2.8 - 2.10) and which may be of independent interest for those working on the inverse Galois problem. For example, it turns out that solvability of (relative) embedding problems over **Q** is a diophantine property (Corollary 5.5).

**Acknowledgement:** This paper was initiated by an invitation of Moshe Jarden to the Tel Aviv University in February/March 2000. I am very greatful to him and to Dan Haran for stimulating discussions on Theorem 1 of the present paper.

# 1 Locally split embedding problems

Let us extract from [Po1], assertion 1.1, and from Erhsov's analysis of 'projective $\Delta^\star$-groups' what seem to be the 'correct' notions for dealing with [strongly] relatively projective profinite groups:

**Definition 1.1** *Let $G$ be a profinite group and let $G_1, \ldots, G_n$ be subgroups of $G$. Then a **locally split embedding problem** for $G$ w.r.t. $G_1, \ldots, G_n$ is given by a pair of epimorphisms $\alpha : G \twoheadrightarrow B$, $\beta : A \twoheadrightarrow B$, where $A$ and $B$ are profinite groups, and by homomorphisms $\beta_i : \alpha(G_i) \to A$ with $\beta \circ \beta_i = id_{\alpha(G_i)}$ $(i = 1, \ldots, n)$. The embedding problem is called **finite**, if $A$ is finite. It is called **reduced**, if $A = <im\,\beta_1, \ldots, im\,\beta_n>$.*

*A **solution** of such a locally split embedding problem is a homomorphism $\gamma : G \to A$ with $\alpha = \beta \circ \gamma$. If $\gamma$ is surjective, it is called a **proper** solution. A solution $\gamma$ is called **locally exact** if $\gamma \mid G_i = \beta_i \circ \alpha \mid G_i$ for each $i = 1, \ldots, n$. And a solution $\gamma$ is called a **locally conjugate** solution if $\gamma(G_i)$ is conjugate to $im\,\beta_i$ in $A$ for each $i = 1, \ldots, n$.*

## 1.1 Characterizing free products via locally split embedding problems

Recall that a profinite group $G$ is the **free product** $G = G_1 \star \cdots \star G_n$ of the subgroups $G_1, \ldots, G_n \leq G$, if any given homomorphisms $\gamma_i : G_i \to H$ $(i = 1, \ldots, n)$ into a profinite group $H$ uniquely extend to a homomorphism $\gamma : G \to H$ (i.e. $\gamma \mid G_i = \gamma_i$ for each $i$). Therefore, every finite locally split embedding problem for $G = G_1 \star \cdots \star G_n$ w.r.t. $G_1, \ldots, G_n$ has a locally exact solution: take $H = A$ and $\gamma_i = \beta_i \circ \alpha \mid G_i$. We now prove the converse:

**Proposition 1.2** *Let $G = <G_1, \ldots, G_n>$ be a profinite group generated by subgroups $G_1, \ldots, G_n \leq G$ and assume that every finite locally split embedding problem for $G$ w.r.t. $G_1, \ldots, G_n$ has a locally exact solution. Then $G = G_1 \star \cdots \star G_n$.*

*Proof:* Let $G'_1, \ldots, G'_n$ be pairwise disjoint isomorphic copies of $G_1, \ldots, G_n$ respectively, fix isomorphisms $\rho_i : G_i \to G'_i$ $(i = 1, \ldots, n)$, let $G' := G'_1 \star \cdots \star G'_n$ and let $\rho : G' \to G$ be the (unique) homomorphism with $\rho \mid G'_i = \rho_i^{-1}$ for $i = 1, \ldots, n$. We first prove the following

**Claim:** *Given any epimorphism $\alpha : G \twoheadrightarrow B$ onto a finite group $B$ and given isomorphisms $\phi_i : \alpha(G_i) \to B'_i$ $(i = 1, \ldots, n)$, where the $B'_i$ are pairwise disjoint, there is a unique epimorphism $\psi : G \twoheadrightarrow B' := B'_1 \star \cdots \star B'_n$ such that for each $i$, $\psi \mid G_i = \phi_i \circ \alpha \mid G_i$.*

To see this, let $\phi : B' \twoheadrightarrow B$ be the unique homomorphism with $\phi \mid B'_i = \phi_i^{-1}$. Then for any open normal subgroup $N \triangleleft B'$ with $N \leq \ker \phi$ we obtain a locally split embedding problem given by $\alpha : G \twoheadrightarrow B$, $\beta : A \twoheadrightarrow B$ and $\beta_i : \alpha(G_i) \to A$, where $A := A_N := B'/N$ with canonical projection $\pi_N : B' \twoheadrightarrow A$ and where $\beta_i := \pi_N \circ \phi_i$ $(i = 1, \ldots, n)$.

By assumption, any such embedding problem has a locally exact solution, so there is some $\gamma = \gamma_N : G \to A$ with $\gamma \mid G_i = \beta_i \circ \alpha \mid G_i$ $(i = 1, \ldots, n)$. Since $B'$ is the inverse limit of these $A_N$, and since the corresponding $\gamma_N$ are unique ($G = <G_1, \ldots, G_n>$) and compatible, the inverse limit $\psi = \lim_{\leftarrow} \gamma_N : G \to B'$ exists and for each $i$, $\psi \mid G_i = \lim_{\leftarrow} \beta_i \circ \alpha \mid G_i = \phi_i \circ \alpha \mid G_i$. As $G = <G_1, \ldots, G_n>$, $\psi$ is also onto and unique, and the claim is proved.

Now associate to any open normal subgroup $M \triangleleft G$ the canonical projection $\alpha_M : G \twoheadrightarrow B_M := G/M$ and the projection $\alpha'_M : G' \twoheadrightarrow B'_M := B'_{M,1} \star \cdots \star B'_{M,n}$ induced by the canonical projections $G'_i \twoheadrightarrow B'_{M,i} := G'_i/\rho_i(G_i \cap M)$ $(i = 1, \ldots, n)$. Then, for each $i$, there is a (unique) isomorphism $\phi_{M,i} : \alpha_M(G_i) \to B'_{M,i}$ such that $\alpha'_M \mid G'_i = \phi_{M,i} \circ \alpha_M \circ \rho_i^{-1}$. The claim thus gives a unique epimorphism $\psi_M : G \twoheadrightarrow B'_M$ with $\psi_M \mid G_i = \phi_{M,i} \circ \alpha_M \mid G_i = \alpha'_M \circ \rho_i$.

Since $G = \lim_{\leftarrow} \alpha_M(G)$, we have $G_i = \lim_{\leftarrow} \alpha_M(G_i)$, so $G'_i = \lim_{\leftarrow} \alpha'_M(G'_i)$ and $G' = \lim_{\leftarrow} \alpha'_M(G'_i) = \lim_{\leftarrow} \alpha'_M(G')$. By the uniqueness of the $\psi_M$ we thus obtain an inverse system of commutative diagrams

$$\begin{array}{ccc} G' & \xrightarrow{\rho} & G \\ \alpha'_M \downarrow & \swarrow \psi_M & \downarrow \alpha_M \\ B'_M & \to & B_M \end{array}$$

with $\psi_M \mid G_i = \alpha'_M \circ \rho_i$. Thus, the inverse limit $\psi := \lim_{\leftarrow} \psi_M : G \to G'$ exists and $\psi \mid G_i = \lim_{\leftarrow} \alpha'_M \circ \rho_i = \rho_i$ for each $i$. Therefore, $\rho \circ \psi = id_G$ and $\rho$ is an isomorphism. $\square$

For the realization of free products of absolute Galois groups as absolute Galois groups (Theorem 1), the free-product criterion of Proposition 1.2 is good enough. In some situations, however, we require a criterion which is easier to check.

**Proposition 1.3** *Let $G = <G_1, \ldots, G_n>$ be a profinite group generated by subgroups $G_1, \ldots, G_n \leq G$. Then $G = G_1 \star \cdots \star G_n$ iff every finite reduced locally split embedding problem for $G$ w.r.t. $G_1, \ldots, G_n$ has a proper locally conjugate solution.*

*Proof:* For the non-trivial direction of the proof it suffices, by Proposition 1.2 that any finite locally split embedding problem for $G$ w.r.t. $G_1, \ldots, G_n$ has a locally exact solution. So let $A, B$ be finite groups, $\alpha : G \twoheadrightarrow B$, $\beta : A \twoheadrightarrow B$ epimorphisms and $\beta_i : \alpha(G_i) \to A$ homomorphisms with $\beta \circ \beta_i = id_{\alpha(G_i)}$ $(i = 1, \ldots, n)$. We have to find a homomorphism $\gamma : G \to A$ with $\gamma \mid G_i = \beta_i \circ \alpha \mid G_i$. Since $B = <\alpha(G_1), \ldots, \alpha(G_n)>)$, we may assume that our embedding problem is reduced.

Let $A_1, \ldots, A_n$ be pairwise disjoint isomorphic copies of $im\, \beta_1, \ldots, im\, \beta_n$ respectively, fix isomorphisms $\pi_i : im\, \beta_i \to A_i$, let $A' = A_1 \star \cdots \star A_n$ and let $\pi : A' \twoheadrightarrow A$ be the epimorphism with $\pi \mid A_i = \pi_i^{-1}$ $(i = 1, \ldots, n)$. Then for any open normal subgroup $N \triangleleft A'$ with $N \leq \ker \pi$ we can canonically lift our given locally split embedding problem from $A$ to $A_N := A'/N$ by setting $\alpha_N = \alpha : G \twoheadrightarrow B$, $\beta_N = \beta \circ \pi^N : A_N \twoheadrightarrow B$ and $\beta_{N,i} = \pi_N \circ \pi_i \circ \beta_i : \alpha(G_i) \to A_N$ $(i = 1, \ldots, n)$, where $\pi^N : A_N \twoheadrightarrow A$ and $\pi_N : A' \twoheadrightarrow A_N$ are the canonical projections (so $\pi = \pi^N \circ \pi_N$):

$$
\begin{array}{ccccccc}
 & & & & & & G \\
 & & & & & & \downarrow \alpha \\
A' & \stackrel{\pi_N}{\to} & A_N & \stackrel{\pi^N}{\to} & A & \stackrel{\beta}{\to} & B \\
\cup & & & & \cup & & \cup \\
A_i & & \stackrel{\pi_i}{\leftarrow} & & im\, \beta_i & \stackrel{\beta_i}{\leftarrow} & \alpha(G_i)
\end{array}
$$

Note that

$$\beta_N \circ \beta_{N,i} = \beta \circ \pi^N \circ \pi_N \circ \pi_i \circ \beta_i = \beta \circ \pi \circ \pi_i \circ \beta_i = \beta \circ \beta_i = id_{\alpha(G_i)},$$

since $\pi \circ \pi_i \mid im\, \beta_i = id_{im\, \beta_i}$,
and that $<im\, \beta_{N,1}, \ldots, im\, \beta_{N,n}> = <\pi_N(A_1), \ldots, \pi_N(A_n)> = A_N$.

By assumption, each of these reduced lifted locally split embedding problems has a proper locally conjugate solution, i.e. there is an epimorphism $\gamma_N : G \twoheadrightarrow A_N$ with $\alpha_N = \beta_N \circ \gamma_N$ such that $\gamma_N(G_i)$ is conjugate to $im\,\beta_{N,i} = \pi_N(A_i)$ for each $i$. Moreover, there are only finitely many locally conjugate solutions $\gamma_N$, because $\gamma_N$ is uniquely determined by the images $\gamma_N(G_i)$ ($i = 1,\ldots,n$): $\gamma_N/G_i = a_i^{-1}(im\,\beta_{N,i})a_i$ for some $a_i \in A_N$ and $\beta_{N,i} \circ \beta_N \mid im\,\beta_{N,i} = id_{im\,\beta_{N,i}}$, so for $g \in G_i$

$$a_i \gamma_N(g) a_i^{-1} = (\beta_{N,i} \circ \beta_N)(a_i \gamma_N(g) a_i^{-1}) = \beta_{N,i}(\beta_N(a_i)\alpha_N(g)\beta_N(a_i)^{-1}).$$

Writing $A' = \lim_{\leftarrow} A_N$ with $N$ ranging over all open normal subgroups of $A'$ with $N \leq \ker \pi$, and writing $\beta' := \beta \circ \pi$, one therefore obtains an epimorphism $\gamma' : G \twoheadrightarrow A'$ with $\alpha = \beta' \circ \gamma'$ such that, for each $i$, $\gamma'(G_i)$ is conjugate to $A_i$ in $A'$. Moreover, for each $i$, $\ker \gamma' \cap G_i = \ker \alpha \cap G_i$, since $\beta'$ is injective on $A_i$ and hence on the conjugate $\gamma'(G_i)$ of $A_i$, and $\ker \alpha \cap G_i = \ker(\pi \circ \beta_i \circ \alpha) \cap G_i$. Therefore, there is an isomorphism $\rho_i : A_i \to \gamma'(G_i)$ with $\gamma' \mid G_i = (\rho_i \circ \pi_i \circ \beta_i \circ \alpha) \mid G_i$ ($i = 1,\ldots,n$).

Now let $\rho : A' \to A'$ be the endomorphism of $A'$ with $\rho \mid A_i = \rho_i$. Surjectivity then passes from $\gamma'$ to $\rho$, and, since $A'$ is finitely generated (hence small), this implies that $\rho$ is an isomorphism.

Thus $\rho^{-1} \circ \gamma' : G \to A'$ is a homomorphism with $\rho^{-1} \circ \gamma' \mid G_i = \pi_i \circ \beta_i \circ \alpha \mid G_i$ ($i = 1,\ldots,n$), and so the induced homomorphism $\gamma = \pi \circ \rho^{-1} \circ \gamma' : G \to A$ is the desired locally exact solution of our locally split embedding problem: for each $i$, $\gamma \mid G_i = \beta_i \circ \alpha \mid G_i$.$\square$

## 1.2 Characterizing relatively projective groups

**Proposition 1.4** *Let $G$ be a profinite group with subgroups $G_1,\ldots,G_n$. Then the following are equivalent:*
*(1) $G$ is [strongly] projective relative to $G_1,\ldots,G_n$.*
*(2) Every embedding problem $\alpha : G \twoheadrightarrow B$, $\beta : A \twoheadrightarrow B$ with a 'local' solution $\gamma_i : G_i \to A$ (i.e. $\alpha \mid G_i = \beta \circ \gamma_i$) for each $i$ has a 'global' solution $\gamma : G \to A$ (i.e. $\alpha = \beta \circ \gamma$) [and $\gamma(G_i)$ is conjugate to $\gamma_i(G_i)$ in $A$ for each $i$].*
*(3) Every locally split embedding problem for $G$ w.r.t. $G_1,\ldots,G_n$ has a [locally conjugate] solution.*
*(4) Every finite locally split embedding problem for $G$ w.r.t. $G_1,\ldots,G_n$ has a [locally conjugate] solution.*

*Moreover, any of these conditions implies [is equivalent to]*
*(5) There is a free profinite group $F$ with $rk\,F \leq rk\,G$ such that there is an*

embedding $\rho : G \hookrightarrow G_1 \star \cdots \star G_n \star F$ [and, for each $i$, $\rho(G_i)$ is conjugate to the factor $G_i$ in $G_1 \star \cdots \star G_n \star F$].

*Proof:* **(1)** $\Rightarrow$ **(2):** Assume **(1)** and let $\alpha : G \twoheadrightarrow B$, $\beta : A \twoheadrightarrow B$ define an embedding problem for $G$ with local solutions $\gamma_i : G_i \to A$ [such that $\gamma(G_i)$ is conjugate to $\gamma_i(G_i)$] for each $i$. Consider the fibre product $A \times_B G := \{(a,g) \in A \times G \mid \beta(a) = \alpha(g)\}$ with the canonical projections $\pi_A : A \times_B G \twoheadrightarrow A$ and $\pi_G : A \times_B G \twoheadrightarrow G$ onto the corresponding coordinate, so $\alpha \circ \pi_G = \beta \circ \pi_A$. Then $\pi_G$ has local splittings

$$\rho_i \quad G_i \to A \times_B G$$
$$g \mapsto (\gamma_i(g), g)$$

for each $i$, and hence, by **(1)**, a global splitting $\rho : G \to A \times_B G$ [with $\rho(G_i)$ conjugate to $\rho_i(G_i)$ in $A \times_B G$, say $\rho(G_i) = (\rho_i(G_i))^{(a_i, g_i)}$]. Now $\gamma := \pi_A \circ \rho$ is the solution we look for: $\alpha = \alpha \circ \pi_G \circ \rho = \beta \circ \pi_A \circ \rho = \beta \circ \gamma$ [and $\gamma(G_i) = \pi_A(\rho(G_i)) = \pi_A((\rho_i(G_i))^{(a_i,g_i)}) = \gamma_i(G_i)^{a_i}$] for each $i$.

**(2)** $\Rightarrow$ **(3):** Just observe that any local splitting $\beta_i : \alpha(G_i) \to A$ for a locally split embedding problem $\alpha : G \twoheadrightarrow B$, $\beta : A \twoheadrightarrow B$ provides a local solution $\gamma_i = \beta_i \circ \alpha \mid G_i$.

**(3)** $\Rightarrow$ **(1):** Any epimorphism $\pi : H \twoheadrightarrow G$ with local splittings $\rho_i : G_i \to H$ defines a locally split embedding problem $\alpha = id_G$, $\beta = \pi$ and $\beta_i = \rho_i$.

**(3)** $\Rightarrow$ **(4):** Clear.

**(4)** $\Rightarrow$ **(1):** Assume **(4)** and let $\pi : H \twoheadrightarrow G$ be an epimorphism with local splittings $\rho_i : G_i \to H$ $(i = 1, \ldots, n)$.

We first proceed as the proof of [Gr], Proposition 1, and consider the case where $\ker \pi$ is finite. For each $i$, $\ker \pi \cap im\, \rho_i = 1$ since $\pi \circ \rho_i = id_{G_i}$. As $\ker \pi$ is finite, there is an open normal subgroup $N \triangleleft H$ with $\ker \pi \cap N im\, \rho_i = 1$ for all $i = 1, \ldots, n$. Hence $N \ker \pi \cap N im\, \rho_i = N$ and so the induced finite embedding problem $\alpha_N : G \twoheadrightarrow G/\pi(N)$, $(\beta_N =) \pi_N : H/N \twoheadrightarrow G/\pi(N)$ where $\alpha_N$ is the canonical projection and $\pi_N(hN) = \pi(h)\pi(N)$ for $h \in H$, is locally split:

$$\ker \pi_N \cap N im\, \rho_i / N = (N \ker \pi \cap N im\, \rho_i)/N = N/N = 1.$$

Thus, there is a solution $\gamma_N : G \to H/N$ with $\pi_N \circ \gamma_N = \alpha_N$.

Now consider the fibre product

$$H/N \times_{G/\pi(N)} G := \{(hN, g) \in H/N \times G \mid \pi_N(hN) = \alpha_N(g)\}$$

and observe that
$$\epsilon: \begin{array}{rcl} H & \to & H/N \times_{G/\pi(N)} G \\ h & \mapsto & (hN, \pi(h)) \end{array}$$
is an embedding, since $\ker \pi \cap N = 1$. Now
$$\rho': \begin{array}{rcl} G & \to & H/N \times_{G/\pi(N)} G \\ g & \mapsto & (\gamma_N(g), g) \end{array}$$
is well-defined ($\pi_N(\gamma_N(g)) = \alpha_N(g)$) and $im\, \rho' \subseteq im\, \epsilon$: for $g \in G$ let $\gamma_N(g) = hN$ for some $h \in H$ and observe that $h$ can be chosen with $\pi(h) = g$ because $\pi(h) \cdot \pi(N) = g \cdot \pi(N)$. Hence $\rho := \epsilon^{-1} \circ \rho'$ is a splitting for $\pi$.

[Moreover, if $\gamma_N$ is a locally conjugate solution, say $\gamma_N(G_i) = (Nim\rho_i(G_i)/N)^{h_i N}$ for some $h_i \in H$, then $\rho'(G_i) = (\epsilon(Nim\, \rho_i(G_i)/N)^{h_i N}$ and so $\rho(G_i)$ is conjugate to $im\, \rho_i$ in $H$.]

Further, let us observe that each [locally conjugate] splitting $\rho$ of $\pi$ can be considered as inverse limit of [l.c.] splittings $\rho_N$ of $\pi_N$, where $N$ runs through all open normal subgroups of $H$ with $\ker \pi \cap Nim\rho = 1$ (again, such $N$'s exist because $\ker \pi$ is finite and $\ker \pi \cap \rho(G) = 1$). Conversely, any such inverse limit of [l.c.] compatible splittings of $\pi_N$ gives a splitting of $\pi$. Since each $\pi_N$ has only finitely many splittings the set of [l.c.] splittings of $\pi$ is an inverse limit of finite sets, and hence compact.

Now let $\ker \pi$ be arbitrary and consider the family $\mathcal{K}$ of normal subgroups $K \triangleleft H$ which are open subgroups of $\ker \pi$. Then for each $K \in \mathcal{K}$, $\pi$ induces a projection $\pi_K : H/K \twoheadrightarrow G$ which splits locally and has a finite kernel, so the set $R_K$ of [locally conjugate] splittings of $\pi_K$ is non-empty and compact. Now the $R_K$ ($K \in \mathcal{K}$) form an inverse system of non-empty compact sets. Hence the inverse limit is non-empty, and any element in it defines a [locally conjugate] splitting of $\pi$

**(1) $\Rightarrow$ (5):** Assuming **(1)**, we can choose a free profinite group $F$ of $rk\, F = rk\, G$ with an epimorphism $\pi_0 : F \twoheadrightarrow G$ extending to an epimorphism $\pi : G_1 \star \cdots \star G_n \star F \to G$ which maps each free factor $G_i$ identically onto the subgroup $G_i$ of $G$. By **(1)**, $\pi$ splits globally [in a locally conjugate way] and any such splitting gives the desired embedding.

**[5] $\Rightarrow$ [1]:** Assume [5]. Then, by the subgroup theorem of Haran ([Ha], Theorem 5.1), $\rho(G)$ is strongly projective relative to $\rho(G_1), \ldots, \rho(G_n)$, since $G_1 \star \cdots \star G_n \star F$ is strongly projective relative to $G_1, \ldots, G_n$. (Compare the following remark). $\square$

**Remark 1.5** *As the proposition shows, our notions of relative [strong] projectivity coincide with Pop's notions (defined via finite embedding problems). Our notion of strong relative projectivity coincides with Haran's notion of relative projectivity and with Ershov's notion of a projective $\Delta$-group: if $G$ is strongly projective relative to subgroups $G_1, \ldots, G_n$ then $G_1, \ldots, G_n$ are separated, i.e. $G_i \cap G_j = 1$ for $i \neq j$. This follows from **(2)**: considering the homomorphism $G_1 \star \cdots \star G_n \star F \to G_1 \times \cdots \times G_n \times F$ that identifies the factors it is clear that for $i \neq j$ any conjugates of $G_i$ and $G_j$ in $G_1 \star \cdots \star G_n \star F$ intersect trivially. This also answers the question implicit in [Er2], Remark 1.*

## 2 Tools from valuation theory

In this section we describe valuation theoretic tools used to realize or to recognize absolute Galois groups as free or projective products of decomposition groups, at the same time introducing (mostly standard) notation and terminology as well as collecting other (mostly well-known) facts. [En] and [Ri] are classical references on valuation theory, the most comprehensive recent book is [K].

### 2.1 Absolutely defectless fields

For a valued field $(F, v)$ we denote valuation ring, maximal ideal, residue field and value group by $\mathcal{O}_v$, $\mathcal{M}_v$, $Fv := \mathcal{O}_v/\mathcal{M}_v$ and $\Gamma_v = v(F^\times) \cong F^\times/\mathcal{O}_v^\times$ respectively. Let $D_v$ be a decomposition subgroup of $G_F$ w.r.t. $v$, i.e. $D_v = G_{F^v}$ for some henselisation $F^v$ of $(F, v)$ in $F^{sep}$. Denoting ramification and inertia subgroup of $D_v$ by $R_v$ and $I_v$, we recall the following

**Facts 2.1** *(a) $R_v$ and $I_v$ are normal subgroups of $D_v$ with $R_v \leq I_v$, and both $R_v$ and $I_v$ have complements in $D_v$ ([KPR]).*
*(b) $D_v/I_v \cong G_{Fv}$*
*(c) $R_v = 1$ if char $Fv = 0$; otherwise (if char $Fv = p > 0$), $R_v$ is a Sylow-$p$ subgroup of $I_v$.*
*(d) The fixed field of $R_v$ in $F^{sep}$ is the smallest subextension of $F^{sep}/F^v$ with separably closed residue field and $q$-divisible value group for all primes $q \neq$ char $Fv$.*

**Definition 2.2** *We call a valued field $(F, v)$* **absolutely defectless** *if $R_v = 1$.*

This definition does not depend on the choice of $D_v$, since any two decomposition subgroups of $G_F$ are conjugate in $G_F$ and this conjugation induces an isomorphism of the corresponding ramification subgroups. Note also, that all finite separable extensions of an absolutely defectless valued field are defectless (i.e. the fundamental equality '$\sum e_i \cdot f_i = n$' holds), but the converse may be false (e.g. $\mathbf{Q}_p$). Moreover, an absolutely defectless field may well have inseparable extensions with defect. (Using the terminology of [K], chapter II.12, $(F, v)$ is absolutely defectless iff the henselisation of $(F, v)$ is 'separably tame'.)

Recall that an extension $(F', v')/(F, v)$ of valued fields is called **immediate** if the canonical embeddings $Fv \hookrightarrow F'v'$ and $\Gamma_v \hookrightarrow \Gamma_{v'}$ are isomorphisms. An extension $(F', v'_1, \ldots, v'_n)/(F, v_1, \ldots, v_n)$ of $n$-fold valued fields is called immediate, if each $(F', v'_i)/(F, v_i)$ is immediate.

**Observation 2.3** *Let $(F', v')/(F, v)$ be an immediate extension of valued fields with $(F', v')$ henselian and absolutely defectless and with $F$ relatively algebraically closed in $F'$. Then $(F, v)$ is also henselian and absolutely defectless and $res : G_{F'} \to G_F$ is an isomorphism.*

*Proof:* Since $F$ is relatively algebraically closed in $F'$, it is clear that $(F, v)$ is henselian and $res : G_{F'} \to G_F$ is surjective. If $F_R$ denotes the fixed field of the ramification subgroup $R_v$ of $G_F$ in $F^{sep}$, then, since $(F', v')/(F, v)$ is immediate, $F'F_R$ is an algebraic extension of $(F', v')$ with separably closed residue field and $q$-divisible value group for all primes $q \neq char\, Fv$. Since $(F', v')$ is absolutely defectless, this implies that $F'F_R = F'^{sep}$ (Fact 2.1(d)), so $F'^{sep} = F'F^{sep}$ and $res : G_{F'} \to G_F$ is injective. Thus, $res$ is an isomorphism and, in particular, $R_v = G_{F_R} \cong G_{F'F_R} = 1$, i.e. $(F, v)$ is also absolutely defectless.□

That $(F, v)$ in the observation is henselian and absolutely defectless was already observed in [K], II., Lemma 12.29 under the weaker assumption that $F'v'/Fv$ be algebraic (instead of $(F', v')/(F, v)$ being immediate).

**Caveat** Note that it may happen that $(F', v')$ is henselian and absolutely defectless, but that a relatively algebraically closed subfield $(F, v)$ is *not* absolutely defectless.

## 2.2 Independence

Two valuations $v$ and $w$ on a field $F$ are called **independent** if they are non-trivial and $F = \mathcal{O}_v \mathcal{O}_w$, i.e., as a ring, $F$ is generated by the proper subrings $\mathcal{O}_v$ and $\mathcal{O}_w$. An important consequence is the well-known

**Approximation Theorem**

*Let $v_1, \ldots, v_n$ be (pairwise) independent valuations on a field $F$. Then, given any $a_1, \ldots, a_n \in F$ and $b_1, \ldots, b_n \in F^\times$, there is an element $x \in F$ with $v_i(x - a_i) > v_i(b_i)$ for all $i = 1, \ldots n$.*

We shall need the following almost trivial

**Observation 2.4** *Two valuations $v$ and $w$ on a field $F$ are independent iff $\forall \gamma \in \Gamma_v\, \exists x \in F : v(x) > \gamma\ \&\ w(x) < 0$.*

*Proof:* If $v$ and $w$ are independent then they are non-trivial. So, given $\gamma \in \Gamma_v$, there are $y, z \in F$ with $v(y) > \gamma$ and $w(z) < 0$. The approximation theorem now provides an element $x \in F$ with $v(x) = v(y) > \gamma$ and $w(x) = w(z) < 0$.

If, conversely, given any $y \in F$, we find an element $x \in F$ with $v(x) > v(y^{-1})$ and $w(x) < 0$, we see that $y = (yx)x^{-1} \in \mathcal{O}_v \mathcal{O}_w$, and that $v$ and $w$ are non-trivial. $\square$

**Corollary 2.5** *Let $v_1, \ldots, v_n$ be independent valuations on $F$ and let $(F', v_1', \ldots, v_n')/(F, v_1, \ldots, v_n)$ be an extension of n-fold valued fields where each $\Gamma_{v_i}$ is cofinal in $\Gamma_{v_i'}$, i.e. $\forall \gamma \in \Gamma_{v_i'}\, \exists \delta \in \Gamma_{v_i}$ with $\delta \geq \gamma$ (e.g. if the extension is immediate or algebraic). Then $v_1', \ldots, v_n'$ are independent.*

For the convenience of the reader, let us reproduce the proof of the following

**Fact 2.6** ([H], Theorem 1.1) *If $F = F_1 \cap \ldots \cap F_n$ for henselian algebraic extensions $F_1, \ldots, F_n$ of a field $F$ inducing independent valuations $v_1, \ldots, v_n$ on $F$, then each $F_i$ is a henselisation of $(F, v_i)$.*

*Proof:* Let $F^{(i)}$ be a henselisation of $(F, v_i)$ in $F_i$ and pick any $\alpha \in F_i$. Then we can approximate the irreducible polynomial of $\alpha$ over $F$ w.r.t. $v_i$ and, for $j \neq i$ some polynomial of the same degree splitting in distinct linear factors over $F$ w.r.t. $v_j$ sufficiently well by some $f \in F[X]$ to guarantee that all zeros of $f$ lie in $\bigcap_{j \neq i} F_j$ and that for some zero $\beta$ of $f$ (close to $\alpha$ in $F_i$) $F^{(i)}(\beta) = F^{(i)}(\alpha) \subseteq F_i$ (Krasner's Lemma). But then $\beta \in \bigcap_j F_j = F$ and $\alpha \in F^{(i)}(\alpha) = F^{(i)}(\beta) = F^{(i)}$. Hence $F_i = F^{(i)}$. $\square$

## 2.3 Decomposition in finite Galois extensions

Recall the following well-known details from ramification theory:

**Facts 2.7** *Let $(F, v)$ be a valued field, let $L/F$ be a finite Galois extension and let $K/F$ be a subextension of $L/F$.*
*(i) $K$ is a decomposition subfield of $L/F$ w.r.t $v$ iff $K = F \cap H$ for some henselisation $H$ of $(F, v)$.*
*(ii) If $K$ is a decomposition subfield of $L/F$ w.r.t. $v$, say $D := Gal(L/K) = \{\sigma \in Gal(L/F) \mid \sigma \mathcal{O}_w = \mathcal{O}_w\}$ for some prolongation $w$ of $v$ to $L$, then $w$ is the only prolongation of $w \mid K$ to $L$ and $v$ has exactly $r := [K : F] = [Gal(L/F) : D]$ distinct prolongations to $L$: the 'conjugates' $w \circ \tau_1, \ldots, w \circ \tau_r$ of $w$, where $Gal(L/F) = \bigcup_{i=1}^{r} \tau_i D$.*
*(iii) If $E/F$ is any finite extension with distinct prolongations $w_1, \ldots, w_k$ of $v$ to $E$, then, given $x_1 \in \mathcal{O}_{w_1}, \ldots, x_k \in \mathcal{O}_{w_k}$, there is some $x \in E$ with $x - x_i \in \mathcal{M}_{w_i}$ for each $i$.*

**Lemma 2.8** *Given a finite Galois extension $L_0/F_0$ and an embedding $\phi : Gal(L_0/F_0) \hookrightarrow A$ into a finite group $A$, there is a Galois extension $(L, w)/(F, v)$ of valued fields with $Lw \cong L_0$, $Fv \cong F_0$, and $Gal(L/F) = A$ such that $\phi(Gal(L_0/F_0))$ is the decomposition subgroup of $A$ w.r.t. $w$. $((L, w)/(F, v)$ is then **tame unramified**, i.e. the inertia subgroup of $Gal(L/F)$ w.r.t. $w$ is trivial.)*

*Proof:* Let $L = F_0(T_1, \ldots, T_d)$, where $d := \sharp A$, say $A = \{1 = a_1, a_2, \ldots, a_d\}$ with $\phi(Gal(L_0/F_0)) = \{a_1, \ldots, a_f\}$ (so $f \mid d$), and where $T_1, \ldots T_d$ are indeterminates over $F_0$.

Let $A$ act on $L$ by acting trivially on $F_0$, and via 'left multiplication' on $\{T_1, \ldots, T_d\}$, i.e. $a_k(T_l) = T_m \Leftrightarrow a_k \cdot a_l = a_m$. Denoting the fixed field under this action of $A$ on $L$ by $F$, it is clear that $L/F$ is a Galois extension with $Gal(L/F) = A$.

Now choose $x \in L_0$ such that $\{\sigma(x) \mid \sigma \in Gal(L_0/F_0)\}$ is a normal base for $L_0/F_0$. Define a ring homomorphism

$$\pi_0 : \begin{array}{rcl} F_0[T_1, \ldots, T_d] & \to & L_0 \\ T_i & \mapsto & \begin{cases} \phi^{-1}(a_i)(x) & \text{for } i \leq f \\ 0 & \text{for } i > f \end{cases} \end{array}$$

with $\pi_0 \mid F_0 = id_{F_0}$ and extend it to a place $\pi : L \to L_0 \cup \{\infty\}$ (Chevalley-Lang) with corresponding valuation $w$. Then $Lw \cong L_0$ and $Kw \cong Fv \cong F_0$, where $K$ is the fixed field of $\phi(Gal(L_0/F_0))$ in $L$ and $v := w \mid F$. Moreover, for $i > f$, $a_i(\mathcal{O}_w) \not\subseteq \mathcal{O}_w$ (e.g. $T_1^{-1} \in \mathcal{O}_w$, but $a_i(T_1^{-1}) = T_i \notin \mathcal{O}_w$). Hence $\phi(Gal(L_0/F_0))$ is the decomposition subgroup of $A$ w.r.t. $w$.□

The lemma, of course, implies that any finite group $A$ with any subgroup $D \leq A$ can be realized as Galois group of a Galois extension of valued fields, where $D$ becomes a decomposition subgroup. And it is not difficult to see that this generalizes to profinite groups. It may, however, be worth noting that this has no analogue for *absolute* Galois groups: in general, not any subgroup $D$ of an absolute Galois group $A = G_F$ can become a decomposition subgroup of $A$ when $A$ is suitably realized as absolute Galois group of some valued field: e.g. if $2 < [A : D] < \infty$ this is not possible.

**Lemma 2.9** *Let $(L, w)/(F, v)$ be a tame unramified finite Galois extension of valued fields, where $Fv$ is infinite. Then there is a primitive element $x \in L = F(x)$ over $F$ with irreducible polynomial $f \in \mathcal{O}_v[X]$ over $F$ such that $\overline{x} \in Lw$ is a primitive element for $Lw/Fv$ and such that $\overline{f} \in Fv[X]$ is the product of the irreducible polynomials of pairwise non-conjugate primitive elements for the Galois extension $Lw/Fv$ over $Fv$. In particular, $f'(x) \in \mathcal{O}_w^\times$.*

*Proof:* Let $K$ be the decompositin subfield of $L/F$ w.r.t. $w$, let $r := [K : F]$, let $\tau_1 = 1, \tau_2, \ldots, \tau_r \in Gal(L/F)$ be representatives for the cosets of $Gal(L/K)$ in $Gal(L/F)$ and let $w_1 := w \circ \tau_1 = w, w_2 := w \circ \tau_2, \ldots, w_r := w \circ \tau_r$ be the distinct prolongations of $v$ to $L$.

Choose $x_1 \in \mathcal{O}_w$ such that $\overline{x_1} \in Lw$ is a primitive element for $Lw = Fv(\overline{x_1})$ over $Fv$, and choose $a_1 := 1, a_2, \ldots, a_r \in \mathcal{O}_v^\times$ with $\overline{a_i} \neq \overline{a_j}$ for $i \neq j$ (this is possible as $Fv$ is infinite). Then $\overline{a_1 x_1}, \ldots \overline{a_r x_1}$ are non-conjugate primitive elements for $Lw/Fv$.

Now choose $x \in L$ with $x - \tau_i(a_i x_1) \in \mathcal{M}_{w_i}$ for $i = 1, \ldots, r$: use Fact 2.7(iii). Since $\overline{\tau_1^{-1}(x)} = \overline{a_1 x_1}, \ldots, \overline{\tau_r^{-1}(x)} = \overline{a_r x_r}$ are non-conjugate primitive elements for $Lw/Fv$, the conjugates $\sigma \tau_i^{-1}(x)$ of $x$ in $L$ over $F$ with $\sigma \in Gal(L/K) \cong Gal(Lw/Fv)$ and $i = 1, \ldots, r$ are all distinct, so $L = L(x)$ (as $[L : F] = r \cdot [L : K]$). Moreover, the irreducible polynomial $f$ of $x$ over $F$ is in $\mathcal{O}_v[X]$, and $\overline{f}$ decomposes over $Fv$ into the product of the irreducible polynomials of $\overline{a_i x_1}$ ($i = 1, \ldots, r$) over $Fv$, which are pairwise coprime, whence $\overline{f'(x)} \neq 0 \in Fv$.□

**Lemma 2.10** *Let $(F, v)$ be a non-trivially valued field, let $L/F$ be a finite Galois extension, let $K/F$ be a subextension of $L/F$ of degree $[K : F] = r$, and assume that $v$ has at most $r$ prolongations to $L$.*

*Then $K$ is a decomposition subfield of $L/F$ iff there is a polynomial $h(X) = X^r + h_{r-1}X^{r-1} + \ldots + h_0 \in F[X]$ which is irreducible over $F$, has a root in $K$ and coefficients with $h_0, \ldots h_{r-2}, 1 + h_{r-1} \in \mathcal{M}_v$.*

*Proof:* '$\Leftarrow$': If $h \in F[X]$ has all the properties mentioned, then $K = F(x)$ for some root $x$ of $h$ in $K$, and, by Hensel's Lemma, $x$ is in some henselisation of $F$ w.r.t. $v$ ($v(h(1)) > 0 = v'(h(1))$), i.e. $K$ is contained in a decomposition subfield $E$ of $L/F$ w.r.t $v$ (Fact 2.7(i)). By Fact 2.7(ii), $v$ has $[E : F]$ prolongations to $L$, so $r \geq [E : F] \geq [K : F] = r$ and hence $K = E$.

'$\Rightarrow$': Let $K$ be a decomposition subfield of $L/K$, say w.r.t. the prolongation $w$ of $v$ to $L$, and let $Gal(L/F) = \bigcup_{i=1}^{r} \tau_i Gal(L/K)$ with $\tau_1 = 1$. Then, by Fact 2.7(ii) and (iii), there is some $y \in K$ with $y \in (1+\mathcal{M}_w) \cap \bigcap_{i=2}^{r} \mathcal{M}_{w \circ \tau_i}$: note that $w \circ \tau_i \mid K \neq w \circ \tau_j \mid K$ for $i \neq j$, since $w$ is the only prolongation of $w \mid K$ to $L$, and so $w \mid K \neq w \circ \tau_j \tau_i^{-1} \mid K$.

Now let $z$ be a primitive element for the extension $K/F(y)$, so $K = F(y, z)$. Since $v$ is non-trivial and $\Gamma_v$ is cofinal in all $\Gamma_{w \circ \tau_i}$, we find infinitely many $a \in \mathcal{M}_v \subseteq F$ such that $az \in \bigcap_{i=1}^{r} \mathcal{M}_{w \circ \tau_i}$. As $K/F$ is finite, there are two such $a$'s, say $a_1 \neq a_2 \in F$, with $F(a_1z+y) = F(a_2z+y)$, so $y, z \in F(a_1z+y)$, and hence $x := a_1z + y$ is a primitive element for $K = F(y, z) = F(x)$ over $F$. Moreover, $x \in (1 + \mathcal{M}_w) \cap \bigcap_{i=2}^{r} \mathcal{M}_{w \circ \tau_i}$, so $\tau_2(x), \ldots, \tau_r(x) \in \mathcal{M}_w$, and the irreducible polynomial of $x$ over $F$:

$$h(X) = \prod_{i=1}^{r} (X - \tau_i(x)) \in F[X]$$

has coefficients $h_0, \ldots, h_{r-2} \in \mathcal{M}_v$ and $h_{r-1} \in -1 + \mathcal{M}_v$.$\square$

## 3 Free products of decomposition groups

The goal of this section is to prove the following proposition and a variation of it (Proposition 3.7).

**Proposition 3.1** *Let $F$ be a field with absolutely defectless valuations $v_1, \ldots, v_n$ and assume that for each $i$, $(F, v_i)$ admits an immediate extension $(\hat{F}_i, \hat{v}_i)/(F, v_i)$* with $\sharp \hat{F}_i > \sharp F$. Let $D_1, \ldots, D_n$ be decomposition subgroups of $G_F$ w.r.t. $v_1, \ldots, v_n$ respectively.

Then there is an immediate extension $(F', v'_1, \ldots, v'_n)/(F, v_1 \ldots, v_n)$ with $\sharp F' = \sharp F$ such that $G_{F'} = D'_1 \star \cdots \star D'_n$, where for each $i$, $D'_i$ is a decomposition subgroup of $G_{F'}$ w.r.t. $v'_i$ and $res : D'_i \to D_i$ is an isomorphism.

**Remark 3.2** *In the hypothesis of the above proposition, we may even assume that for each $i$, $(\hat{F}_i, \hat{v}_i)$ is henselian, is absolutely defectless and contains the henselisation $F_i$ of $(F, v_i)$ corresponding to $D_i = G_{F_i}$. Moreover, we may assume that $G_F = <D_1, \ldots, D_n>$, i.e. $F = F_1 \cap \ldots \cap F_n$.*

*Proof of the remark:* First pass to a henselisation $(\hat{F}'_i, \hat{v}'_i)$ of $(\hat{F}_i, \hat{v}_i)$. This is an immediate extension and $\hat{F}'_i \cap F^{sep}$ is a henselisation of $(F, v_i)$: since $(F, v_i)$ is absolutely defectless any immediate separably algebraic henselian extension is a henselisation. So there is an isomorphism $\hat{F}'_i \cap F^{sep} \to F_i$ over $(F, v_i)$ which extends to an isomorphism $(\hat{F}'_i, \hat{v}'_i) \to (\hat{F}''_i, \hat{v}''_i)$ over $(F, v_i)$. Thus, $(\hat{F}''_i, \hat{v}''_i)/(F, v_i)$ is an immediate henselian extension containing $F_i$.

Now we pass to the fixed field $\hat{F}'''_i$ of a complement of the ramification subgroup of $G_{\hat{F}''_i}$ (cf. Fact 2.1(a)) to obtain an absolutely defectless field. $(\hat{F}'''_i, \hat{v}'''_i)/(F, v_i)$ is still immediate since passing from $\hat{F}''_i$ to $\hat{F}'''_i$ means passing to the perfect hull for the residue field and passing to the $p$-divisible hull for the value group ($p = char\, Fv$, cf. Fact 2.1(d)). As $(F, v_i)$ is absolutely defectless, $Fv \cong \hat{F}''_i$ is already perfect and $\Gamma_{v_i} \cong \Gamma_{\hat{v}''_i}$ $p$-divisible.

Finally, replacing $F$ by $F_1 \cap \ldots \cap F_n$ doesn't change any of the hypotheses: note that the assumption that $(\hat{F}_i, \hat{v}_i)/(F, v_i)$ be an immediate extension of higher cardinality implies that the valuations $v_1, \ldots, v_n$ are all non-trivial, so all fields involved are infinite and hence, all algebraic extensions of the same cardinality.$\square$

Before proving the proposition let us first isolate the key arguments in three lemmas.

**Lemma 3.3** *Given an n-fold valued field $(F, v_1, \ldots, v_n)$ and decomposition groups $D_1, \ldots, D_n$ as in Proposition 3.1 and the remark thereafter, there is an immediate extension $(F', v'_1, \ldots, v'_n)/(F, v_1 \ldots, v_n)$ with $\sharp F' = \sharp F$ such that $v'_1, \ldots, v'_n$ are independent, absolutely defectless, and $G_{F'} = <D'_1 \ldots, D'_n>$ for decomposition subgroups $D'_i$ of $G_{F'}$ w.r.t. $v'_i$ for which $res : D'_i \to D_i$ is an isomorphism ($i = 1, \ldots, n$).*

*Proof:* Let $X = X_1 \dot\cup \ldots \dot\cup X_n$ be a (partitioned) set of indeterminates over $F$ with $\sharp X_1 = \cdots = \sharp X_n = \sharp F$. For each $i$, $(\hat{F}_i, \hat{v}_i)/(F, v_i)$ is immediate and $\sharp \hat{F}_i > \sharp F$, so we find an embedding $\phi_i : F(X) \to \hat{F}_i$ such that $\hat{v}_i(\phi_i(x)) < 0$ for all $x \in X_j$ ($j \neq i$) and such that $\forall \gamma \in \Gamma_{v_i} \exists x \in X_i$ with $\hat{v}_i(\phi_i(x)) > \gamma$.

After passing to isomorphic copies of $(\hat{F}_i, \hat{v}_i)$ over $(F_i, v_i)$ (as in the proof of the previous remark), we may for all $i = 1, \ldots, n$, $x \in X$ identify $\phi_i(x)$ with $x$, so that then $X \subseteq \bigcap_i \hat{F}_i$, and for each $i$, $\hat{v}_i(X_j) < 0$ when $j \neq i$ and $\hat{v}_i(X_i)$ is cofinal in $\Gamma_{\hat{v}_i}$. Hence, by Observation 2.4, the valuations on $F(X)$ induced by $\hat{v}_1, \ldots, \hat{v}_n$ are independent.

Passing, if necessary, once more to isomorphic copies of $(\hat{F}_i, \hat{v}_i)$ over $(F_i(X), \hat{v}_i \mid F_i(X))$, we may even assume that the relative algebraic closures $F'_i$ of $F_i(X)$ in $\hat{F}_i$ are all contained in a fixed algebraic closure of $F(X)$. By Observation 2.3, each $F'_i$ is then henselian and absolutely defectless (w.r.t. $\hat{v}_i \mid F'_i$).

Now let $F' = F'_1 \cap \ldots \cap F'_n$ and, for each $i$, let $v'_i := \hat{v}_i \mid F'$. Then, by Corollary 2.5, $v'_1, \ldots v'_n$ are independent, and so, by Fact 2.6, each $F'_i$ is a henselisation of $(F', v'_i)$. Hence, by Observation 2.3, $res : D'_i := G_{F'_i} \to D_i$ is an isomorphism for each $i$: $F_i$ is relatively algebraically closed in $\hat{F}_i$, and so in $F'_i$.□

**Lemma 3.4** *([HJ], Proof of Prop. 14.1, Part C, [P1], Thm. 3.1, and [Er2], Lemma 1) Let $L/K$ be a finite Galois extension with Galois group $B = Gal(L/K)$ and let $\beta : A \twoheadrightarrow B$ be an epimorphism of finite groups. Consider the elements of $A$ as indeterminates over $L$ and let $A$ act on $L(A)$ via group multiplication on $A$ and via the given Galois action of $\beta(A) = B$ on $L$. Assume that $D \leq A$ is a subgroup with $D \cap \ker \beta = 1$.*

*Then $L(A)^D := \{x \in L(A) \mid \sigma x = x \text{ for all } \sigma \in D\}$ is purely transcendental over $L^{\beta(D)}$.*

**Lemma 3.5** *Let $(F, v_1, \ldots, v_n)$ and $D_1, \ldots, D_n$ satisfy the assumptions of Proposition 3.1 and the remark thereafter. Assume, in addition, that either $v_1, \ldots, v_n$ are independent or that char $Fv_i = 0$ for all $i$. Let $\alpha : G_F \twoheadrightarrow B$, $\beta : A \twoheadrightarrow B$ and $\beta_i : \alpha(D_i) \to A$ ($i = 1, \ldots, n$) be the data for a locally split embedding problem for $G_F$ w.r.t. $D_1, \ldots D_n$, i.e. $A$, $B$ are finite groups, $\alpha$, $\beta$ are epimorphisms, $\beta_1, \ldots, \beta_n$ are homomorphisms with $A = <im \beta_1, \ldots, im \beta_n>$ and, for each $i$, $\beta \circ \beta_i = id_{\alpha(D_i)}$.*

*Then $(F, v_1, \ldots, v_n)$ admits an immediate extension $(F', v'_1, \ldots, v'_n)$ with $\sharp F' = \sharp F$, where each $(F', v'_i)$ is absolutely defectless and where $G_{F'}$ has* decomposition subgroups $D'_1, \ldots, D'_n$ w.r.t. $v'_1, \ldots, v'_n$ respectively such that $G_{F'} =< D'_1, \ldots D'_n >$, $res : D'_i \to D_i$ is an isomorphism for each $i$ and such that the lifted locally split embedding problem for $G_{F'}$ w.r.t. $D'_1, \ldots, D'_n$ (i.e. $\alpha' := \alpha \circ res : G_{F'} \twoheadrightarrow B$, $\beta' := \beta : A \to B$, $\beta'_i := \beta_i \circ res : \alpha'(D'_i) \to A$) has a locally exact solution (i.e. $\exists \gamma' : G_{F'} \twoheadrightarrow A$ with $\gamma' \mid D'_i = \beta'_i \circ \alpha' \mid D'_i$).

*Proof:* Let $L = (F^{sep})^{\ker \alpha}$ be the fixed field of $\ker \alpha$ in $F^{sep}$, so $L/F$ is a Galois extension with $Gal(L/F) \cong B$ and (identifying those two groups) we may assume that $\alpha = res : G_F \to Gal(L/F) = B$. Then, for each $i$, $L^{\alpha(D_i)} = F_i \cap L$.

Consider the field $L(A)$ with the $A$-action from the previous lemma. Then $L(A)/L(A)^A$ is a Galois extension with Galois group $A$ and, for each $i$, the conclusion of the lemma (with $D = im\,\beta_i$) says that the extension $L(A)^{im\,\beta_i}/L^{\alpha(D_i)}$ is purely transcendental.

Since $L^{\alpha(D_i)} \subseteq \hat{F}_i$ and $\sharp \hat{F}_i > \sharp F = \sharp L^{\alpha(D_i)}$, we can consider $L(A)^{im\,\beta_i}$ as subfield of $\hat{F}_i$.

As before, passing, if necessary, to isomorphic copies of $\hat{F}_i$ over $F_i L(A)^{im\,\beta_i}$, we may even assume that all relative algebraic closures $F'_i$ of $F_i L(A)^{im\,\beta_i}$ in $\hat{F}_i$ are contained in a fixed algebraic closure of $L(A)$. Then, by Observation 2.3, each $F'_i$ is henselian and absolutely defectless w.r.t. $\hat{v}_i \mid F'_i$, $res : D'_i := G_{F'} \to D_i$ is an isomorphism, and $(F' := F'_1 \cap \ldots \cap F'_n, v'_1, \ldots, v'_n)$ with $v'_i := \hat{v}_i \mid F'$ is an immediate extension of $(F, v_1, \ldots, v_n)$ with $G_{F'} =< D'_1, \ldots, D'_n >$ and $\sharp F' = \sharp F$.

Moreover, each $F'_i$ is a henselisation of $(F', v'_i)$ (so the $D'_i$ are decomposition groups): if $v_1, \ldots, v_n$ are independent, then so are, by Corollary 2.5, $v'_1, \ldots, v'_n$, and hence, by Fact 2.6, $F'_i$ must be a henselisation of $(F', v'_i)$; and if $char\,Fv_i = 0$ for all $i$, then as immediate henselian algebraic extension, $F'_i$ is, again, a henselisation of $(F'_i, v'_i)$.

Now $\gamma' := res : G_{F'} \to Gal(L(A)/L(A)^A) = A$ is a homomorphism with

$$\gamma' \mid D'_i = res : D'_i \twoheadrightarrow Gal(L(A)/(F'_i \cap L(A))) = im\,\beta_i,$$

so

$$\beta_i^{-1} \circ \gamma' \mid D'_i = \alpha' \mid D'_i = res : D'_i \twoheadrightarrow Gal(L/(F_i \cap L)) = \alpha(D_i) \leq B.$$

Hence, $\gamma'$ is a locally exact solution of the lifted locally split embedding problem.$\square$

*Proof of Proposition 3.1:* The proof of the proposition is now a standard chain construction. We may assume from the start that $(F, v_1 \ldots, v_n)$ etc. satisfies the conditions in the remark following the proposition. By Lemma 3.3, we may also assume that $v_1, \ldots, v_n$ are independent (this is not necessary if *char $Fv_i = 0$* for all $i$).

We first find an immediate absolutely defectless extension $(F^1, v_1^1, \ldots, v_n^1)/(F, v_1, \ldots, v_n)$ with $\sharp F^1 = \sharp F$ and with decomposition subgroups $D_1^1, \ldots, D_n^1$ of $G_{F^1}$ w.r.t. $v_1^1 \ldots, v_n^1$ respectively such that $G_{F^1} = <D_1^1, \ldots, D_n^1>$, such that $res : D_i^1 \to D_i$ is an isomorphism for each $i$, and such that each locally split embedding problem for $G_F$ w.r.t. $D_1, \ldots, D_n$ has a locally exact solution when lifted to $G_{F^1}$ w.r.t. $D_1^1, \ldots, D_n^1$.

This is achieved by an ordinal enumeration of these embedding problems $(EP_\kappa)_{\kappa < \lambda}$ and constructing an (ordinal) chain of absolutely defectless immediate extensions $(F_\kappa, v_{\kappa,1}, \ldots, v_{\kappa,n})/(F, v_1 \ldots, v_n)$ with $\sharp F_\kappa = \sharp F$, with decomposition subgroups $D_{\kappa,i}$ and isomorphisms $res : D_{\kappa,i} \to D_i$ ($i = 1, \ldots, n$) such that all $EP_\mu$ for $\mu < \kappa$ have a solution when lifted to $G_{F_\kappa}$. For successor ordinals this is done by Lemma 3.5, and for limit ordinals by taking unions of the fields constructed 'before' ($\sharp F_\kappa$ never increases since there are only $\sharp F$-many locally split embedding problems for $G_F$). Note that for $\mu < \nu < \kappa$, solutions for $EP_\mu$ lift from $G_{F_\nu}$ to $G_{F_\kappa}$. Then $(F^1, v_1^1, \ldots, v_n^1) := (F_\lambda, v_{\lambda,1} \ldots, v_{\lambda,n})$ has all the required properties.

Now we iterate this process and construct $F^2, F^3, \ldots$ (with valuations, decomposition groups, same cardinality etc.) solving all locally split embedding problems for $G_{F^i}$ in $G_{F^{i+1}}$ and let $F' = \bigcup_{i=1}^\infty F^i$ (with $v_1', \ldots$ etc.). Then $\sharp F' = \sharp F$, $G_{F'} = <D_1', \ldots, D_n'>$ and each locally split embedding problem for $G_{F'}$ w.r.t. $D_1', \ldots, D_n'$ has a locally exact solution. Hence, by Proposition 1.2, $G_{F'} = D_1' \star \cdots \star D_n'$.□

As a consequence, we shall now prove a variant of Proposition 3.1, dropping details about immediacy and absolute defect both from hypothesis and conclusion, but retaining the Galois theoretic data. The reduction of Proposition 3.7 to Proposition 3.1 proceeds via the following

**Lemma 3.6** *Given a field $K$, there is an absolutely defectless henselian valued field $(L, w)$ with $\sharp L = \max\{\sharp K, \aleph_0\}$ admitting an immediate extension $(\hat{L}, \hat{w})$ with $\sharp \hat{L} > \sharp L$ such that $K \subseteq \mathcal{O}_w$, $Lw$ is the perfect hull of $K$, $\Gamma_w$ is divisible, and, hence, $res : G_L \to G_K$ is an isomorphism.*

*Proof:* Choose an infinite set $X$ of indeterminates over $K$ with $\sharp X = \max\{\sharp K, \aleph_0\}$ and fix some well-ordering '$<$' on $X$. Let $L = K(X)$ and let $w$ be the '$(X, <)$-**adic**' valuation on $L$: For any finite subset $\{x_1 < k_2 < \ldots < x_n\} \subseteq X$, the restriction of $w$ to $K(x_1, \ldots x_n)$ is (equivalent to) the composed valuation $w_{x_1} \oplus \cdots \oplus w_{x_n}$, where $w_{x_i}$ is the $x_i$-adic valuation on the rational function field $K(x_{i+1}, \ldots, x_n)(x_i)$ in $x_i$ over $K(x_{i+1}, \ldots, x_n)$. To make this consistent, define
$$\Gamma_w := \bigoplus_{x \in X} \mathbf{Z} \cdot \gamma_x,$$
where '$\bigoplus$' is the lexicographic sum w.r.t. the (well-)ordering induced by $<$ on $X$ under the bijection $(x \mapsto \gamma_x)_{x \in X}$, define $w(x) := \gamma_x$ for all $x \in X$ and define $w$ to be trivial on $K$. This uniquely determines a valuation $w$ on $L$ with residue field $K$ and value group $\Gamma_w$.

Now the field of formal Laurent series in $X$ over $K$
$$\hat{L} := K((X)) := \{\alpha = \sum_{\gamma \in \Gamma_w} a_\gamma t^\gamma \mid a_\gamma \in K \ \& \ supp(\alpha) \text{ is well-ordered}\}$$
(where $supp(\alpha) := \{\gamma \in \Gamma_w \mid a_\gamma \neq 0\}$) with the canonical henselian valuation
$$\hat{w}(\alpha) := \min\{\gamma \mid a_\gamma \neq 0\} \text{ for all } \alpha = \sum_{\gamma \in \Gamma_w} a_\gamma t^\gamma \in \hat{L}$$
is an immediate extension of $(L, w)$. (Henselianity of $\hat{w}$ is proved, e.g., in [PrC] II.5 Satz 4 and III.2 Satz 17.) Moreover, $\sharp \hat{L} > \sharp L$, since any subset $Y \subseteq X$ is well-ordered and, thus, gives an element
$$\alpha_Y := \sum_{x \in Y} x = \sum_{x \in Y} t^{\gamma_x} \in \hat{L},$$
where $\alpha_Y \neq \alpha_{Y'}$ for $Y \neq Y' \subseteq X$. So $\sharp \hat{L} \geq \sharp\{Y \subseteq X\} > \sharp X = \sharp L$.

Finally, replace $(\hat{L}, \hat{w})$ by the fixed field of the complement of the inertia subgroup of $G_{\hat{L}}$, and replace $(L, w)$ by its relative algebraic closure in the new (absolutely defectless henselian) field $(\hat{L}, \hat{w})$. Then, by Observation 2.3, $(L, w)$ is absolutely defectless and henselian, $(\hat{L}, \hat{w})/(L, w)$ is still immediate, now with divisible value group and the perfect hull of $K$ as residue field. Since $K \subseteq \mathcal{O}_w$, there is a commutative diagram

$$\begin{array}{ccc} G_L & \stackrel{\cong}{\Rightarrow} & G_{Lw} \\ res \downarrow & & \downarrow \cong res \\ G_K & = & G_K \end{array}$$

and so $res : G_L \to G_K$ is an isomorphism as well.$\square$

It may be worth noting that it was only by the special choice of $\Gamma_w$ that $\sharp K((\Gamma_w)) > \sharp K + \sharp \Gamma_w$. If $\Gamma_w = \mathbf{R}$, for example, $\sharp K((\mathbf{R})) = \sharp K + \sharp \mathbf{R}$, since well-ordered subsets of $\mathbf{R}$ are countable.

We conclude this section by a variant of Proposition 3.1:

**Proposition 3.7** *Let $(F, v_1, \ldots, v_n)$ be an n-fold valued field and let $D_1, \ldots, D_n$ be decomposition subgroups of $G_F$ w.r.t. $v_1, \ldots, v_n$ respectively. Then there is an extension $(F', v'_1, \ldots, v'_n)/(F, v_1, \ldots, v_n)$ with $\sharp F' = \max\{\sharp F, \aleph_0\}$ such that $G_{F'} \cong D'_1 \star \cdots \star D'_n$, where for each $i$, $D'_i$ is a decomposition subgroup of $G_{F'}$ w.r.t. $v'_i$, and $res : D'_i \to D_i$ is an isomorphism.*

**Remark 3.8** *The extension of n-fold valued fields established in the proposition is, in general, no longer immediate, but (as the proof will show), there is a divisible ordered abelian group $\Gamma$ such that, for each $i$,*

$$\Gamma_{v'_i} = \Gamma \oplus_{lex} p^{-\infty}\Gamma_{v_i},$$

*where $p = \operatorname{char} F$ (not necessarily $\operatorname{char} Fv_i$!) and $p^{-\infty}\Gamma_{v_i}$ denotes the p-divisible hull of $\Gamma_{v_i}$ ($:= \Gamma_{v_i}$ if $p = 0$), and $F'v'_i = Fv_i$ or the perfect hull of $Fv_i$, again depending on whether or not $\operatorname{char} F = 0$.*

*Proof:* Let $(L, w)$ be the field constructed in the previous lemma from $K = F$, and let $w_1 = \ldots = w_n = w$. By Proposition 3.1, there is an immediate extension $(L', w'_1, \ldots, w'_n)/(L, w_1, \ldots, w_n)$ with $\sharp L' = \sharp L = \max\{\sharp F, \aleph_0\}$ such that $G_{L'} = D_{w'_1} \star \cdots \star D_{w'_n}$, where for each $i$, $D_{w'_i}$ is a decomposition subgroup of $G_{L'}$ w.r.t. $w'_i$, and $\rho_i := res : D_{w'_i} \to D_{w_i} = G_L \cong G_F$ is an isomorphism.

Now consider for each $i$, the composed valuation $w'_i \oplus v_i$, i.e. the refinement of $w'_i$ by the unique prolongation of (again denoted by) $v_i$ from $F$ to the perfect hull $L'w'_i \cong Lw_i = Lw$ of $F$. Let $D'_i := \rho_i^{-1}(D_i)$, let $F'_i$ be the fixed field of $D'_i$ and let $F' = F'_1 \cap \ldots \cap F'_n$.

Then $F'_i$ is a henselisation of $(L', w'_i \oplus v_i)$, hence of the intermediate field $F'$ w.r.t. the induced prolongation $v'_i$ of $w'_i \oplus v_i$ from $L'$ to $F'$, and $res : D'_i := G_{F'_i} \to D_i$ is an isomorphism. Clearly, $\sharp F' = \max\{\sharp F, \aleph_0\}$. And since the subgroup of a free product of profinite groups $G_1 \star \cdots \star G_n$ generated by subgroups $H_i \leq G_i$ ($i = 1, \ldots, n$) is the free product of these subgroups: $<H_1, \ldots, H_n> = H_1 \star \cdots \star H_n$, we also have $G_{F'} = D'_1 \star \cdots D'_n$.$\square$

## 4 Proof of Theorem 1

**1.** It clearly suffices to prove Theorem 1 for $n = 2$. So we are given two fields $F_1$, $F_2$ and we want to find a field $F$ with $G_F \cong G_{F_1} \star G_{F_2}$, where $char\ F = char\ F_1$, provided $char\ F_1 = char\ F_2$.

**2.** It is well-known that the absolute Galois group of a field $K$ of characteristic $p > 0$ can be realized as absolute Galois group of a field $L$ of characteristic 0: just make $K$ the residue field of a valuation of mixed characteristic (extend the $p$-adic valuation on $\mathbf{Q}$ canonically to a valuation on the purely transcendental extension $\mathbf{Q}(X)$ of $\mathbf{Q}$ with residue field $\mathbf{F}_p(X)$, where $X$ is a transcendence basis of $K$ over $\mathbf{F}_p$ and adjoin roots of minimal polynomials of all elements of $K$ over $\mathbf{F}_p(X)$ lifted to $\mathbf{Q}(X)$), pass to the henselisation $L'$ of $L$ and then to the fixed field of a complement of the inertia subgroup of $G_{L'}$ (use Fact 2.1(a)). Hence we may assume that $char\ F_1 = char\ F_2$.

**3.** Since for any subgroups $H_1 \leq G_1$ and $H_2 \leq G_2$ of profinite groups $G_1$, $G_2$, the subgroup generated by $H_1$ and $H_2$ (under the canonical embeddings of $G_1$, $G_2$) in $G_1 \star G_2$ is $H_1 \star H_2$, and since subgroups of absolute Galois groups are absolute Galois groups, it suffices to realize $G_{k(X)} \star G_{k(Y)}$ as absolute Galois group of a field of the characteristic of $k$, where $k = \mathbf{Q}$ or $\mathbf{F}_p$ is the prime subfield of $F_1$ and $F_2$, and where $X$ and $Y$ are transcendence bases of $F_1$ and $F_2$ over $k$ which we may assume algebraically independent.

**4.** Now choose valuations $v_1$ and $v_2$ on $F = k(X \cup Y)$ with $Fv_1 = k(X)$ and $Fv_2 = k(Y)$ and apply Proposition 3.7 to obtain a field with absolute Galois group $D_1 \star D_2$, where $D_i$ is a decomposition subgroup of $G_F$ w.r.t. $v_i$ ($i = 1, 2$). Passing to complements of the inertia subgroups of $D_1$, $D_2$ we obtain, once more applying the argument in 3., a field with absolute Galois group $G_{k(X)} \star G_{k(Y)}$.□

By the same arguments as in Step 3 and 4, it is clear that the proof of Theorem 1 can be reduced to realizing $G_{k(X)} \star G_{k(X)}$ as absolute Galois group having the same characteristic as the prime field $k$, where $X$ is in *infinite* set of indeterminates over $k$. But we do not know whether $G_{k(X)} \cong G_{k(X)} \star G_{k(X)}$.

Theorem 1 has an almost trivial generalisation to 'pro-$\mathcal{C}$ Galois groups': Let $\mathcal{C}$ be an **almost full family** of finite groups, i.e. $\mathcal{C}$ is closed under homomorphic images, subgroups and direct products. A **pro-$\mathcal{C}$ group** is then an inverse limit of groups in $\mathcal{C}$ and the **free pro-$\mathcal{C}$ product** $G_1 \star_\mathcal{C} \cdots \star_\mathcal{C} G_n$ of pro-$\mathcal{C}$ groups $G_1, \ldots, G_n$ is a pro-$\mathcal{C}$ group $G$ admitting embeddings $\epsilon_i : G_i \to G$ such that given any homomorphisms $\gamma_i : G_i \to H$ into a pro-$\mathcal{C}$ group $H$ there is a unique homomorphism $\gamma : G \to H$ with $\gamma_i = \gamma \circ \epsilon_i$ $(i = 1, \ldots, n)$.
The **pro-$\mathcal{C}$ Galois group** $G_F(\mathcal{C})$ of $F$ is the maximal pro-$\mathcal{C}$ quotient of $G_F$, i.e. the Galois group of the compositum of all finite Galois extensions of $F$ with Galois group in $\mathcal{C}$.

**Corollary 4.1** *Given any fields $F_1, \ldots, F_n$, there is a field $F$ with $G_F(\mathcal{C}) \cong G_{F_1}(\mathcal{C}) \star_\mathcal{C} \cdots \star_\mathcal{C} G_{F_n}(\mathcal{C})$. Moreover, $F$ can be chosen to have the same characteristic as all $F_i$, provided they have the same characteritic.*

*Proof:* This is, because the maximal pro-$\mathcal{C}$ quotient of the free product is the free pro-**C** product of the maximal pro-$\mathcal{C}$ quotients.□

# 5 Regularly closed fields

## 5.1 Realizing free products over regularly closed fields

**Proposition 5.1** *Let $(F, v_1, \ldots, v_n)$ be an $n$-fold valued field with corresponding henselisations $F_1, \ldots, F_n$. Then there is an extension $(F', v'_1, \ldots, v'_n)/(F, v_1, \ldots, v_n)$ of $n$-fold valued fields such that $F'$ is regularly closed w.r.t. $v'_1, \ldots, v'_n$, the $v'_i$ are independent and, for each $i$, there is a henselisation $F'_i$ containing $F_i$ where $res : G_{F'_i} \to G_{F_i}$ is an isomorphism and $G_{F'} \cong G_{F'_1} \star \cdots \star G_{F'_n}$.*

*Proof:* In the proof of Proposition 3.7 we already constructed an extension $(F', v'_1, \ldots, v'_n)/(F, v_1, \ldots, v_n)$ of $n$-fold valued fields and, for each $i$ a henselisation $F'_i$ of $(F', v'_i)$ containing $F_i$ such that $res : G_{F'_i} \to G_{F_i}$ is an isomorphism and $G_{F'} \cong G_{F'_1} \star \cdots \star G_{F'_n}$. Moreover, by the way how Proposition 3.1 was used in the proof of Proposition 3.7, $v'_1, \ldots, v'_n$ were independent and, for each $i$, there was an immediate embedding $F'_i \hookrightarrow \hat{F}_i$ of absolutely defectless henselian valued fields where $\sharp \hat{F}_i > \sharp F'_i$, and, again, $res : G_{\hat{F}_i} \to G_{F'_i}$ was an isomorphism.

It now suffices to find a regularly closed extension $(F'', v''_1, \ldots, v''_n)/(F', v'_1, \ldots, v'_n)$ with $F'' = F''_1 \cap \ldots \cap F''_n$ where each $(F'', v''_i)$ embeds (as valued field) into $\hat{F}_i$ and where $F''_i$ is the henselisation of $(F'', v''_i)$ corresponding to the relative algebraic closure of the embedding of $F''$ in $\hat{F}_i$. For then $res : G_{F''_i} \to G_{F'_i}$ is an isomorphsim for each $i$, and so $res : G_{F''} \to G_{F'}$ is an isomorphism as well: surjectivity is clear since $G_{F'}$ is generated by the $G_{F'_i}$. And since $G_{F'} \cong G_{F'_1} \star \cdots \star G_{F'_n}$ there is a (unique) homomorphism $res^{-1} : G_{F'} \to G_{F''}$ induced by the local inverses $res_i^{-1} : G_{F'_i} \to G_{F''_i}$ with $res \circ res^{-1} = id_{G_{F'}}$. As $G_{F''}$ is generated by the $G_{F''_i}$, $res^{-1}$ is surjective and hence $res$ is injective.

To find $(F'', v''_1, \ldots, v''_n)$ one proceeds exactly as in the proof of [HP], Theorem 3.1. By [HP], Theorem 1.8, it suffices to satisfy the local-global principle for rational points on affine *plane curves*. So using a standard chain argument (as in the proof of Proposition 3.1), the crucial step is to find, given an absolutely irreducible polynomial $f(X,Y) \in F'[X,Y]$ with a simple zero $(a_i, b_i)$ in each $F'_i$, a regular $n$-fold valued field extension of $(F', v'_1, \ldots, v'_n)$ of the cardinality of $F'$ which embeds in each $\hat{F}_i$ (as valued field w.r.t. the corresponding prolongation of $v'_i$) and which has a zero of $f$.

To achieve this, it obviously suffices to embed the function field $F'(x,y)$ (where $(x, y)$ is a generic point of the curve) into each $\hat{F}_i$. But this is easy since $\sharp\hat{F}_i > \sharp F_i$: the point $(a_i, b_i)$ is simple, say $\frac{df}{dY}(a_i, b_i) \neq 0$; choose $\epsilon \in \hat{F}_i \setminus F'_i$ with $\hat{v}_i(\epsilon)$ big enough to guarantee that $\hat{v}_i(f(a_i + \epsilon, b_i)) > 2\hat{v}_i(\frac{df}{dY}(a_i + \epsilon, b_i))$; then $a'_i := a_i + \epsilon \in \hat{F}_i$ is transcendental over $F'$ and, since $\hat{F}_i$ is henselian, there is some $b'_i \in \hat{F}_i$ with $f(a'_i, b'_i) = 0$; the embedding of $F'(x,y)$ over $F'$ into $\hat{F}_i$ is now given by mapping $(x, y) \mapsto (a'_i, b'_i)$.□

**Corollary 5.2** *Let $G$ be a profinite group which is strongly projective relative to subgroups $G_1, \ldots, G_n$ and assume that each $G_i$ is isomorphic to an absolute Galois group. Then there is a regularly closed n-fold valued field $(F, v_1, \ldots, v_n)$ with $v_1, \ldots, v_n$ independent and there is an isomorphism $\phi : G \to G_F$ such that for each $i$, $\phi(G_i)$ is a decomposition subgroup of $G_F$ w.r.t. $v_i$.*

*Proof:* We will first construct an $n$-fold valued field $(F, v_1, \ldots, v_n)$ satisfying all stated properties except being regularly closed. To this end let $G_{n+1}$ be a free profinite group with $rk\, G_{n+1} = rk\, G$. Then $G_{n+1}$ is an absolute Galois group as well and we find, as in the proof of Theorem 1, an $(n+1)$-fold valued field $(K, w_1, \ldots, w_{n+1})$ with $G_K \cong G_1 \star \cdots \star G_{n+1}$, where the free factors $G_i$ are decomposition subgroups of $G_K$ w.r.t. $w_i$, where the $w_i$ are independent and where each $(K, w_i)$ allows immediate extensions of higher cardinality. Let $\pi : G_K \to G$ be an epimorphism identifying the free factors $G_1, \ldots, G_n$ with the correpsonding subgroups of $G$ and projecting $G_{n+1}$ onto $G$. As $G$ is strongly projective relative to $G_1, \ldots G_n$, there is a splitting $\phi : G \to G_K$ of $\pi$ with $\phi(G_i)$ conjugate to the factor $G_i$ in $G_K$ (for $i = 1, \ldots, n$). Now let $(F, v_1, \ldots, v_n)$ be the fixed field of $\phi(G)$ where each $v_i$ is induced from the henselisation $F_i$ of $(K, w_i)$ with $G_{F_i} = \phi(G_i)$.

Now we continue as in the previous proof. Given a curve $\mathcal{C}$ over $F$ with $F_i$-rational points we find immediate prolongations of $v_1, \ldots, v_n$ to the function field $L = F(x, y)$ of $\mathcal{C}$ with henselisations $L_i$ containing $F_i$ such that $res : G_{L_i} \to G_{F_i}$ is an isomorphism for each $i$. By strong projectivity, again, the epimorphsim $res : G_L \to G_F$ has a splitting $\phi : G_F \to G_L$. Let $(F', v'_1, \ldots, v'_n)$ be the fixed field of $\phi(G_F)$ with $v'_i$ induced from henselisations $F'_i$ with $G_{F'_i} = \phi(G_{F_i})$. Then $res : G_{F'} \to G_F$ and $res : G_{F'_i} \to G_{F_i}$ ($i = 1, \ldots, n$) are isomorphisms and $\mathcal{C}$ has an $F'$-rational point.

Finally, again, a chain of such extensions leads to the desired regularly closed field. □

## 5.2 The absolute Galois group of regularly closed fields

### 5.2.1 Solvability of embedding problems is existential

Let us recall that a (finite) **embedding problem** for a profinite group $G$ is given by a pair of epimorphisms $\alpha : G \twoheadrightarrow B$, $\beta : A \twoheadrightarrow B$, (where $A$ and $B$ are finite groups). A **solution** resp. a **proper solution** of the embedding problem is a homomorphism resp. an epimorphism $\gamma : G \to A$ such that $\alpha = \beta \circ \gamma$.

**Observation 5.3** *Solvability of a finite embedding problem for the absolute Galois group $G_F$ of a field $F$ can be expressed by an existential (first-order) formula in the language of fields with parameters from $F$.*

Let $\alpha : G_F \twoheadrightarrow B$, $\beta : A \twoheadrightarrow B$ be the data of a finite embedding problem for $G_F$. Let $E := Fix \ker \alpha$ be the fixed field of the kernel of $\alpha$ and let $\psi : Gal(E/F) \to B$ be the unique isomorphism making the diagram

$$\begin{array}{ccc} G_F & = & G_F \\ res_{F^{sep}/E} \downarrow & & \downarrow \alpha \\ Gal(E/F) & \stackrel{\psi}{\to} & B \end{array}$$

commute. We first prove the following
**Claim:** *The above embdding problem has a proper solution iff there is a*

*Galois extension $L/F$ containing $E$ and an isomorphism $\phi : Gal(L/F) \to A$ such that the following diagram commutes:*

$$\begin{array}{ccc} Gal(L/F) & \xrightarrow{\phi} & A \\ res_{L/E} \downarrow & & \downarrow \beta \\ Gal(E/F) & \xrightarrow{\psi} & B \end{array}$$

To prove the claim, assume first that the embedding problem has a proper solution $\gamma : G_F \twoheadrightarrow A$. Then the field $L := Fix \ker \gamma$ is a Galois extension of $F$ containing $E$ (since $\ker \gamma \subseteq \ker \alpha$) and there is a unique isomorphism $\phi : Gal(L/F) \to A$ making the top square of the following diagram commute:

$$\begin{array}{ccc} G_F & = & G_F \\ res_{F^{sep}/L} \downarrow & & \downarrow \gamma \\ Gal(L/F) & \xrightarrow{\phi} & A \\ res_{L/E} \downarrow & & \downarrow \beta \\ Gal(E/F) & \xrightarrow{\psi} & B \end{array}$$

Since $\alpha = \beta \circ \gamma$ and $res_{F^{sep}/E} = res_{L/E} \circ res_{F^{sep}/L}$, the outer square commutes as well and, hence, so does the bottom square.

For the converse, assume there is a Galois extension $L/F$ containing $E$ and an isomorphism $\phi : Gal(L/F) \to A$ making the bottom square of the diagram above commute. We define $\gamma := \phi \circ res_{F^{sep}/L}$, so that the top square commutes. Then $\gamma : G_F \to A$ is an epimorphism and

$$\alpha = \psi \circ res_{F^{sep}/E} = \psi \circ res_{L/E} \circ res_{F^{sep}/L} = (\beta \circ \phi) \circ (\phi^{-1} \circ \gamma) = \beta \circ \gamma.$$

The claim is proved.

Our next step is to express the existence of a Galois extension $L/F$ as described in the claim by an 'almost existential' formula in the language of fields $\{+, \times, 0, 1\}$ (allowing parameters from $F$). We first express that there is a Galois extension $L/F$ with $Gal(L/F) \cong A$. This is equivalent to the existence of an irreducible (monic) polynomial $f \in F[X]$ of degree $d := \sharp A$ such that the $F$-algebra $F[X]/(f)$ contains $d$ distinct zeros $x_1, \ldots, x_d$ of $f$, where for each $k$, the map $x_1 \mapsto x_k$ induces the permutation $\sigma_k$ of $\{x_1, \ldots, x_d\}$ corresponding to the permutation of $A = \{a_1 = 1, a_2, \ldots, a_d\}$ given by left multiplication with $a_k$: then $L = F(x_1) = \cdots = F(x_n)$ is a Galois extension with $Gal(L/F) \cong A$. The elements of the $F$-algebra $F[X]/(f)$ can be regarded as $d$-tuples of elements of $F$

$$r_0 + r_1 X + \ldots + r_{d-1} X^{d-1}/(f) \mapsto (r_0, \ldots, r_{d-1}),$$

addition is componentwise and multiplication is expressible via polynomials in the coefficients (only depending on the coefficients of $f$). Therefore, the existence of a Galois extension $L/F$ with $Gal(L/F) \cong A$ is equivalent to the formula

$$\exists c, x_1, \ldots, x_d, u \in F^d : \Phi(c, x_1, \ldots, x_d, u), \text{ where}$$

$$\begin{aligned}\Phi(c, x_1, \ldots, x_d, u) = & \bigwedge_{k=1}^d f_c(x_k) = 0 \wedge u \cdot \prod_{k \neq l}(x_k - x_l) = 1 \\ & \wedge \bigwedge_{k,l,m: a_k \cdot a_l = a_m} x_{l,0} + x_{l,1} x_k + \ldots + x_{l,d-1} x_k^{d-1} = x_m \\ & \wedge f_c \text{ is irreducible.}\end{aligned}$$

Here $x_l = (x_{l,0}, \ldots, x_{l,d-1})$, and for $c = (c_0, \ldots, c_{d-1}) \in F^d$, '$f_c$' denotes the polynomial $f_c(X) = X^d + c_{d-1} X^{d-1} + \ldots + c_0$. Of course, addition and multiplication on $F^d$ occuring in the formula is induced from the $F$-algebra $F[X]/(f_c)$ and $1 = (1, 0, \ldots, 0) \in F^d$.

To express all properties of $L/F$ required by the claim, let $E = F(\zeta)$ and let $g \in F[X]$ be the irrducible polynomial of $\zeta$ over $F$, say $deg\ g = e$, and let $\zeta_1 = \zeta, \zeta_2, \ldots, \zeta_e$ be the conjugates of $\zeta$ over $F$. Then the existence of a Galois extension $L/F$ containing $E$ and of an isomorphism $\phi : Gal(L/F) \to A$ with $\beta \circ \phi = \psi \circ res_{L/E}$ is equivalent to the formula

$$\exists c, x_1, \ldots, x_d, u, z_1 \ldots, z_e \in F^d : \Phi(c, x_1, \ldots, x_d, u) \wedge \Psi(c, x_1, \ldots, x_d, z_1, \ldots, z_e),$$

where

$$\begin{aligned}\Psi(c, x_1, \ldots, x_d, z_1, \ldots, z_e) = & \bigwedge_{j=1}^e g(z_j) = 0 \wedge \bigwedge_{j \neq j'} z_j \neq z_{j'} \\ & \wedge \bigwedge_{j,j',j'': \psi^{-1}(\beta(a_{j''}))(\zeta_{j'})=\zeta_j} z_{j',0} + z_{j',1} x_{j''} + \ldots + z_{j',d-1} x_{j''}^{d-1} = z_j\end{aligned}.$$

Again, addition and multiplication on $F^d$ occuring in the formula are inherited from $F[X]/(f_c)$, and $z_{j'} = (z_{j',0}, \ldots, z_{j',d-1})$.

The formula which expresses proper solvability of our embedding problem is existential except for the phrase '$f_c$ is irreducible' in the formula $\Phi$. It may now come as a minor surprise that the most naive way of making the formula existential — delete '$f_c$ is irreducible' — works, provided $F$ is infinite. Yet for finite fields $F$, $G_F \cong \hat{\mathbf{Z}}$ is projective, and so every embedding problem has a solution (and any formula true in $F$ is equivalent to that truth). Hence, from now on, $F$ is assumed to be infinite.

So let us first assume that the new formula holds for $F$, say with $c, x_1, \ldots, x_d$, $u, z_1, \ldots z_e \in F^d$ witnessing this. Then $A$ may be considered as group of $F$-algebra automorphisms of the $d$-dimensional $F$-algebra

$$L' := F[X]/(f_c) = F(x_1) = \cdots = F(x_d).$$

For each irreducible factor $f$ of $f_c$ over $F$, the canonical $F$-algebra epimorphism

$$\pi : L' \twoheadrightarrow L := F[X]/(f)$$

maps the $d$ distinct zeros $x_1, \ldots, x_d$ of $f_c$ in $L'$ to $d$ distinct zeros of $f_c$ in the field $L$:

$$\pi(u) \cdot \prod_{k \neq l}(\pi(x_k) - \pi(x_l)) = \pi(u \cdot \prod_{k \neq l}(x_k - x_l)) = \pi(1) = 1.$$

Hence, $f_c$ is a separable polynomial and

$$L = F(\pi(x_1)) = \cdots = F(\pi(x_d))$$

is a Galois extension of $F$: the splitting field of $f_c$ over $F$. (So all irreducible factors of $f_c$ over $F$ generate the *same* field extension $L/F$.) Moreover, $E \subseteq L$, since $g(\pi(z_1)) = \pi(g(z_1)) = 0$. It is also easily checked that

$$G := \{a \in A \mid a(ker\pi) \subseteq ker\pi\} = \{a \in A \mid f(\pi(a(X + (f_c)))) = 0\} \leq A$$

and that

$$\pi^\star : G \to Gal(L/F)$$
$$a \mapsto \sigma_a : \begin{array}{ccc} L & \to & L \\ X + (f) & \mapsto & \pi(a(X + (f_c))) \end{array}$$

is a well-defined isomorphism: note that $A$ acts simply transitive on the zeros of $f_c$ in $L'$, that $G$ acts simply transitive on the zeros of $f_c$ in $L'$ which become zeros of $f$ in $L$, and so $\pi^\star(G)$ acts simply transitive on the zeros of $f$ in $L$. So $\phi := (\pi^\star)^{-1}$ embeds $Gal(L/F)$ into $A$, and the formula $\Psi$ implies that $\beta \circ \phi = \psi \circ res_{L/E}$.

For the converse, assume $L/F$ is a Galois extension containing $E$ with an embedding $\phi : Gal(L/F) \to A$ such that $\beta \circ \phi = \psi \circ res_{L/E}$. Now use

Lemma 2.8 and 2.9 to find a tame unramified Galois extension $(M, w)/(K, v)$ of valued fields with $Gal(M/K) \cong A$, $Mw \cong L$, $Kv \cong F$ and a primitive element $x \in M = K(x)$ over $K$ with irreducible polynomial $f \in \mathcal{O}_v[X]$ over $K$ such that $f'(x) \in \mathcal{O}_w^\times$. Then the polynomial $f_c := \overline{f} \in F[X]$ is separable and the action of $A$ on the zeros of $f$ in $M$ induces an action of $A$ on the zeros of $f_c$ in the $F$-algebra $F[X]/(f_c)$ fulfilling our formula:

If $w_1, \ldots, w_r$ are the distinct prolongations of $v$ to $M$ then all zeros of $f$ are in $\mathcal{O}_{w_1} \cap \ldots \cap \mathcal{O}_{w_r}$, and $A$ permutes the $d$ distinct zeros of $f$ in each $\mathcal{O}_{w_j}$. Thus $A$ acts simply transitive on the corresponding (tuples of) zeros of $f$ in the $\mathcal{O}_v$-algebra $\mathcal{O}_{w_1} \times \cdots \times \mathcal{O}_{w_r}$, and, via the canonical ring epimorphism

$$\mathcal{O}_{w_1} \times \cdots \times \mathcal{O}_{w_r} \twoheadrightarrow Mw_1 \times \cdots \times Mw_r \cong F[X]/(f_c),$$

on the corresponding zeros $x_1, \ldots, x_d$ of $f_c$ in $F[X]/(f_c)$. In particular, $\prod_{k \neq l}(x_k - x_l)$ is a unit in $F[X]/(f_c)$ since $\prod_{a_k \neq a_l \in A}(a_k(x) - a_l(x)) \in \mathcal{O}_{w_1}^\times \cap \ldots \cap \mathcal{O}_{w_r}^\times$. □

The idea to consider 'Galois-*algebra* extensions' rather than just field extensions when dealing with embedding problems for absolute Galois groups already occurs in Hasse's 1948-paper [Hs], section 1.

**Proposition 5.4** *Let $(F, v_1, \ldots, v_n)$ be an n-fold valued field with henselisations $F_1, \ldots, F_n \subseteq F^{sep}$. Then 'locally conjugate solvability' of a finite locally split embedding problem for $G_F$ w.r.t. $G_{F_1}, \ldots, G_{F_n}$ can be expressed by an existential formula in the language of n-fold valued fields with parameters from $F$.*

*Proof:* Let $\alpha : G_F \twoheadrightarrow B$, $\beta : A \twoheadrightarrow B$, $\beta_i : \alpha(G_{F_i}) \to A$ $(i = 1, \ldots, n)$ be the data of a finite locally split embedding problem for $G_F$ w.r.t. $G_{F_1}, \ldots, G_{F_n}$, i.e. $A$, $B$ are finite groups, $\alpha$, $\beta$ are epimorphisms, and the $\beta_i$ are homomorphisms with $\beta \circ \beta_i = id_{\alpha(G_{F_i})}$.

As in the proof of Observation 5.3, we let $E = Fix \ker \alpha$ be the fixed field of the kernel of $\alpha$, we let $\psi : Gal(E/F) \to B$ be the unique isomorphism making the diagram

$$\begin{array}{ccc} G_F & = & G_F \\ res_{F^{sep}/E} \downarrow & & \downarrow \alpha \\ Gal(E/F) & \stackrel{\psi}{\to} & B \end{array}$$

commute and first prove the following
**Claim:** *The above embedding problem has a proper locally conjugate solution*

*iff there is a Galois extension $L/F$ containing $E$ and an isomorphism $\phi : Gal(L/F) \to A$ such that the diagram*

$$\begin{array}{ccc} Gal(L/F) & \stackrel{\phi}{\to} & A \\ res_{L/E} \downarrow & & \downarrow \beta \\ Gal(E/F) & \stackrel{\psi}{\to} & B \end{array}$$

*commutes and $\phi^{-1}(im\,\beta_i)$ is a decomposition subgroup of $Gal(L/F)$ w.r.t. $v_i$ for each $i$.*

'$\Rightarrow$': Since a proper locally conjugate solution $\gamma$ of a locally split embedding problem is, in particular, a proper solution of the embedding problem given by $\alpha, \beta$, we can construct $L$ and $\phi$ as in the proof of Observation 5.3. So we only have to check that, for each $i$, $\phi^{-1}(im\,\beta_i)$ is a decomposition subgroup of $Gal(L/F)$ w.r.t. $v_i$. But this follows since

$$\phi^{-1}(\gamma(G_{F_i})) = res_{F^{sep}/L}(G_{F_i}) = Gal(L/L \cap F_i)$$

is a decomposition subgroup of $Gal(L/F)$ w.r.t. $v_i$ and is conjugate to $\phi^{-1}(im\,\beta_i)$ in $Gal(L/F)$.

'$\Leftarrow$': For the converse, the claim in the proof of Observation 5.3 already provides an epimorphism $\gamma : G_F \twoheadrightarrow A$ with $\alpha = \beta \circ \gamma$. Then, for each $i$, both $\phi^{-1}(im\,\beta_i)$ and $\phi^{-1}(\gamma(G_{F_i})) = Gal(L/L \cap F_i)$ are decomposition subgroups of $Gal(L/F)$ w.r.t. $v_i$. Hence, they are conjugate in $Gal(L/F)$ and so $im\,\beta_i$ and $\gamma(G_{F_i})$ are conjugate in $A$, i.e. $\gamma$ is a locally conjugate solution and the claim is proved.

Since any solution $\gamma$ of an embedding problem is a *proper* solution of the modified embedding problem where $A$ is replaced by $im\,\gamma$, the claim provides also a criterion for locally conjugate solvability of the given locally split embedding problem: it is the criterion of the claim except that $\phi$ is only required to be an embedding such that, for each $i$, there is some $a_i \in A$ with $(im\,\beta_i)^{a_i} \subseteq im\,\phi$ and with the property that $\phi^{-1}((im\,\beta_i)^{a_i})$ is a decomposition subgroup of $Gal(L/F)$ w.r.t. $v_i$.

Our next step, again as in the proof of Observation 5.3, is to express proper locally conjugate solvability by an 'almost existential' formula in the language of $n$-fold valued fields, and we may take the formula

$$\exists c, x_1, \ldots, x_d, u, z_1, \ldots, z_e \in F^d : \Phi(c, x_1, \ldots, x_d, u) \wedge \psi(c, x_1 \ldots, x_d, z_1, \ldots, z_e)$$

from the proof of Observation 5.3 to express existence of a Galois extension $L/F$ containing $E$ and an isomorphism $\phi : Gal(L/F) \to A$ such that $\psi \circ res_{L/E} = \beta \circ \phi$. So we only have to express that, for each $i$, $\phi^{-1}(im\,\beta_i)$ is a decomposition subgroup of $Gal(L/F)$ w.r.t. $v_i$.

To achieve this for one $i$, let $K$ be the fixed field of $\phi^{-1}(im\,\beta_i)$ in $L$. Then $K \cap E$ is the fixed field of $\psi^{-1}(\alpha(G_{F_i}))$ in $E$, i.e. $K \cap E = F_i \cap E$ is a decomposition subfield of $E/F$. Hence $v$ has exactly $[K \cap E : F]$ prolongations to $E$ (Fact 2.7 (ii)), and thus at most $r := [L : E] \cdot [K \cap E : F]$ prolongations to $L$. Since $\beta_i : \alpha(G_{F_i}) \to A$ is injective, $[L : K] = \sharp im\,\beta_i = \sharp \alpha(G_{F_i}) = [E : K \cap E]$, and so

$$\begin{aligned} r &= [L : E] \cdot [K \cap E : F] = \tfrac{[L:K\cap E]}{[E:K\cap E]} \cdot [K \cap E : F] \\ &= \tfrac{[L:K]\cdot[K:K\cap E]}{[E:K\cap E]} \cdot [K \cap E : F] = [K : F]. \end{aligned}$$

Lemma 2.10 now says that $K$ is a decomposition subfield of $L/F$ (i.e. that $\phi^{-1}(im\,\beta_i)$ is a decomposition subgroup of $Gal(L/F)$) w.r.t. $v_i$ iff there is a polynomial $h(X) = X^r + t_{r-1}X^{r-1} + \ldots + t_0 \in F[X]$ such that $h$ is irreducible over $F$, has a root in $K$ and satisfies $t_0, \ldots, t_{r-2}, 1 + t_{r-1} \in \mathcal{M}_{v_i}$.

Under the identification $L \cong F[X]/(f_c) \cong F^d$ used in the formula already imported from the proof of Observation 5.3, this is expressed, for each $i$, by the existential formula

$$\exists y_i \in F^d \exists h_i, w_i \in F^{r_i} \Theta_i(c, x_1, \ldots, x_d, y_i, h_i, w_i), \text{ where}$$

$$\begin{aligned} \Theta_i(c, x_1, \ldots, x_d, y_i, h_i, w_i) = &\bigwedge\nolimits_{k : a_k \in im\,\beta_i} a_k(y_i) = y_i \\ &\wedge \prod\nolimits_{k \neq l \in S_i} (a_k(y_i) - a_l(y_i)) \in F^\times \\ &\wedge y_i^{r_i} + h_{i,r_i-1} y_i^{r_i-1} + \ldots + h_{i,0} = 0 \\ &\wedge h_{i,0}, \ldots, h_{i,r_i-2}, 1 + h_{i,r_i-1} \in \mathcal{M}_{v_i} \\ &\wedge g_i(w_{i,0} + w_{i,1} y_i + \ldots + w_{i,r_i-1} y_i^{r_i-1}) = 0 \\ &\wedge \bigvee\nolimits_{p,q \in \{1,\ldots,d\}^d} det((z_1^{p_j} y_i^{q_j})_{1 \leq j \leq d}) \neq 0, \end{aligned}$$

where $a_k(y_i) := y_{i,0} + y_{i,1} x_k + \ldots + y_{i,d-1} x_k^{d-1}$, where $S_i \subseteq \{1, \ldots, d\}$ such that the $a_k$ with $k \in S_i$ form a system of coset representatives of $A/im\,\beta_i$, where $r_i := \sharp A / \sharp im\,\beta_i$, and where $g_i$ is the irreducible polynomial of a primitive element of $F_i \cap E$ over $F$. Note that the first two lines of $\Theta_i$ express that $F(y_i)$ is the fixed field of $\phi^{-1}(im\,\beta_i)$ and that hence the polynomial $X^{r_i} + h_{i,r_i-1} X^{r_i-1} + \ldots + h_{i,0}$ is the *irreducible* polynomial of $y_i$ over $F$. Moreover, formulas of the type $t \in \mathcal{O}_v^\times$ and $t \in \mathcal{M}_v$ are existential in the language of valued fields:

$$t \in \mathcal{O}_v^\times \iff t \in \mathcal{O}_v \wedge \exists t' \in \mathcal{O}_v : t \cdot t' = 1$$
$$t \in \mathcal{M}_v \iff t \in \mathcal{O}_v \wedge \exists t' \notin \mathcal{O}_v : t \cdot t' = 1.$$

The existential formula expressing locally conjugate solvability of our locally split embedding problem is now

$$\exists c, x_1, \ldots x_d, u, y_1, \ldots, y_n, z_1, \ldots, z_e \in F^d \exists h_1, w_1 \in F^{r_1} \ldots \exists h_n, w_n \in F^{r_n} :$$

$$\Phi'(c, x_1, \ldots, x_d, u) \wedge \Psi(c; x_1, \ldots, x_d, z_1, \ldots, z_e) \wedge \bigwedge_{i=1}^n \Theta_i(c, x_1, \ldots, x_d, y_i, h_i, w_i),$$

where $\Phi'$ is the formula $\Phi$ except that the phrase '$f_c$ is irreducible' is deleted. Note that we may assume $F$ to be infinite since otherwise all valuations are trivial and hence any locally split embedding problem trivially solvable.

If the existential formula holds, then we obtain a solution of our locally split embedding problem, as in the proof of Observation 5.3, by considering the $F$-algebra $L' := F[X]/(f_c)$, an irreducible factor $f$ of $f_c$ over $F$, and the $F$-algebra epimorphism $\pi : L' \twoheadrightarrow L := F[X]/(f)$ onto the Galois extension $L/F$ with $\mathrm{Gal}(L/F) \overset{\phi}{\hookrightarrow} A \leq \mathrm{Aut}_F L'$, where $\beta \circ \phi = \psi \circ \mathrm{res}_{L/E}$. We have to check that, for each $i$, $\mathrm{im}\,\phi$ contains a conjugate of $\mathrm{im}\,\beta_i$ in $A$, and that the corresponding subgroup of $\mathrm{Gal}(L/F)$ is a decomposition subgroup of $\mathrm{Gal}(L/F)$.

So let us focus on one $i$, and observe that the polynomial $h_i \in F[X]$ of degree $r_i$ has $r_i$ distinct roots in $L'$: the elements $a_k(y_i)$ with $k \in S_i$ and that the condition in the second line of $\Theta_i$ is carried over to $L$ by $\pi$. Hence, $h_i$ has $r_i$ distinct zeros in $L$. It is easily checked that exactly one irreducible factor $h$ of $h_i$ over $F$ is of the shape $X^r + t_{r-1}X^{r-1} + \ldots + t_0$ with $t_0, \ldots, t_{r-2}, 1 + t_{r-1} \in \mathcal{M}_{v_i}$. So we may assume (after replacing $y_i$ by a suitable $a_k(y_i)$ — this is fixed by a subgroup conjugate to $\mathrm{im}\,\beta_i$ in $A$, namely $(\mathrm{im}\,\beta_i)^{a_k}$) that $h(\pi(y_i)) = 0$.

Then $K := F(\pi(y_i))$ contains a zero of $g_i$ (i.e. a decomposition subfield of $E/F$ w.r.t. $v_i$), and is itself (by Hensel's lemma applied to $h$) contained in some henselisation of $(F, v_i)$. Hence $K \cap E$ is a decomposition subfield of $E$ w.r.t. $v_i$. Moreover, the last line of $\Theta_i$ means that $L' = F(z_1, y_i)$, so

$$L = \pi(L') = \pi(F(z_1))\pi(F(y_i)) = EK.$$

Thus, $res : Gal(L/K) \to Gal(E/E \cap K)$ is an isomorphism and $K$ is a decomposition subfield of $L/F$. Finally, $\phi(Gal(L/K)) \leq im\, \beta_i$, as $Gal(L/K)$ fixes $\pi(y_i)$, and since

$$\sharp Gal(L/K) = \sharp Gal(E/E \cap K) = \sharp \alpha(G_{F_i}) = \sharp im\, \beta_i,$$

equality holds.

Conversely, assume that $L/F$ is a Galois extension containing $E$, that $\phi : Gal(L/K) \to A$ is an embedding with $\beta \circ \phi = \psi \circ res_{L/E}$ and that, for each $i$, there is some $a_i \in A$ such that $(im\, \beta_i)^{a_i} \subseteq im\, \phi$ and such that $\phi^{-1}((im\, \beta_i)^{a_i})$ is a decomposition subgroup of $Gal(L/F)$ w.r.t. $v_i$. Once again quoting from the proof of Observation 5.3, we find a tame unramified Galois extension $(M, w)/(K, v)$ of valued fields with $Gal(M/K) \cong A$, $Mw \cong L$, $Kv \cong F$ and a primitive element $x \in M = K(x)$ over $K$ with irreducible polynomial $f \in \mathcal{O}_v[X]$ over $K$ such that $f'(x) \in \mathcal{O}_w^\times$ and such that the action of $A$ on the distinct zeros of the polynomial $f_c := \overline{f} \in F[X]$ in the $F$-algebra $F[X]/(f_c)$ satisfies the formulas $\Phi'$ and $\Psi$.

So only $\Theta_i$ remains to be checked for each $i$. To this end, choose $y \in L$ such that $F(y)$ is the fixed field of $(im\, \beta_i)^{a_i}$ in $L$ and the irreducible polynomial of $y$ over $F$ is of the form $X^r + t_{r-1}X^{r-1} + \ldots + t_0$ with $t_0, \ldots, t_{r-2}, 1 + t_{r-1} \in \mathcal{M}_{v_i}$ (Lemma 2.10). Let $w_1, \ldots, w_s$ be the distinct prolongations of $v$ to $M$, choose distinct $a_2, \ldots, a_s \in \mathcal{M}_{v_i} \setminus \{0\}$, and, by Fact 2.1(iii), choose $\tilde{y}_i \in M$ with

$$\tilde{y}_i \in (y + \mathcal{M}_{w_1}) \cap (a_2 y + \mathcal{M}_{w_2}) \cap \ldots \cap (a_s y + \mathcal{M}_{w_s}).$$

Then $K(\tilde{y}_i)$ is the fixed field of $(im\, \beta_i)^{a_i}$ in $M$ and the irreducible polynomial $\tilde{h}_i \in K[X]$ of $\tilde{y}_i$ over $K$ has coefficients in $\mathcal{O}_v$ and the induced polynomial $h_i := \overline{\tilde{h}_i}$ over $F$ has coefficients $h_{i,0}, \ldots, h_{i,r_i-2}, 1 + h_{i,r_i-1} \in \mathcal{M}_{v_i}$. Moreover, $h_i$ has $r_i$ distinct zeros in $L$ (which gives the second line of $\Theta_i$). And, finally, $L' = F(z_1, y_i)$ since $M = K(\tilde{z}_1, \tilde{y}_i)$, where $\tilde{z}_1 \in \mathcal{O}_w$ is a lifting of $z_1$ in $M$. $\square$

As an application of independent interest let us mention the following Corollary, which might be helpful for the question whether the inverse Galois problem for $\mathbf{Q}$ is decidable: note that there is still a chance that the existential theory of $\mathbf{Q}$ be decidable.

**Corollary 5.5** *Let $F$ be a number field and let $F_1, \ldots, F_n$ be henselisations w.r.t. $n$ distinct primes on $F$. Then locally conjugate solvability of a locally* split embedding problem for $G_F$ w.r.t. $G_{F_1}, \ldots, G_{F_n}$ is a diophantine property, i.e. equivalent to an existential first-order formula in the language of fields.

*Proof:* This is just a combination of our Proposition 5.4 with Rumely's existential definability of valuation rings in number fields ([Ru]).□

Another immediate consequence is the following

**Corollary 5.6** *Let $(F, v_1, \ldots, v_n)$ and $(F', v'_1, \ldots, v'_n)$ be n-fold valued fields which are elementarily equivalent in the the language of n-fold valued fields. Then $G_F$ is the free product of suitable decomposition subgroups (w.r.t. $v_1, \ldots, v_n$) iff the same is true for $G_{F'}$.*

*Proof:* This is immediate from the claim in the proof of the above proposition together with Proposition 1.3.□

### 5.2.2 Proof of Theorem 2

In order to prove Theorem 2 we will use the following model theoretic characterisation of regularly closed fields:

**Fact 5.7** *An n-fold valued field $(F, v_1, \ldots, v_n)$ with pairwise independent valuations and with corresponding henselisations $F_1, \ldots, F_n$ is regularly closed w.r.t. $v_1, \ldots, v_n$ iff $(F, v_1, \ldots, v_n)$ is existentially closed (in the language of n-fold valued fields) in any regular extension $(F', v'_1, \ldots v'_n)$ provided each $F_i$ is existentially closed in some (any) henselisation $F'_i$ of $(F', v'_i)$ containing $F_i$ (in the language of valued fields).*

*Proof:* This is Theorem 1.9 together with Theorem 4.1 of [HP]: note that the general assumption made in [HP] that $F$ should be of characteristic 0 does not enter the proof of those two Theorems.□

*Proof of Theorem 2:* Let $(F, v_1, \ldots, v_n)$ be a regularly closed $n$-fold valued field where $v_1, \ldots, v_n$ are pairwise independent, and let $F_1, \ldots, F_n$ be corresponding henselisations. In order to show that $G_F$ is relative projective w.r.t $G_{F_1}, \ldots, G_{F_n}$, we have to show, by Proposition 1.4, that any finite locally split embedding problem for $G_F$ w.r.t $G_{F_1}, \ldots, G_{F_n}$ has a locally conjugate solution.

So let $\alpha : G_F \twoheadrightarrow B$, $\beta : A \twoheadrightarrow B$, $\beta_i : \alpha(G_{F_i}) \to A$ ($i = 1, \ldots, n$) be the data of such a finite locally split embedding problem. Using the above fact together with Proposition 5.4, it suffices to find an extension $(F', v'_1, \ldots, v'_n)/(F, v_1, \ldots, v_n)$ and, for each $i$, a henselisation $F'_i$ of $(F', v'_i)$ containing $F_i$, such that $F'/F$ is regular, $F_i$ is existentially closed in $F'_i$ (in the language of valued fields) and such that the locally split embedding problem (lifted via $res: G_{F'} \twoheadrightarrow G_F$ to $G_{F'}$) has a locally conjugate solution.

To achieve this, we have to adjust the proof of Lemma 3.5 to our situation. Let, for each $i$, $F_i^\star \succeq F_i$ be an elementary extension (of valued fields) such that $\sharp F_i^\star > \sharp F_i$. Let $L$ be the fixed field of $\ker \alpha$ in $F^{sep}$, and, after the canonical identification, $Gal(L/F) = B$, $\alpha = res: G_F \to Gal(L/F) = B$ and $L^{\alpha(G_{F_i})} = L \cap F_i$.

Consider the field $L' := L(A)$ with the $A$-action of Lemma 3.4, and let $F' := L(A)^A$ be the fixed field under this action. Then $L'/F'$ is a Galois extension with Galois group $A$, $F'/F$ is regular, and, for each $i$, $L'^{im \beta_i}/L \cap F_i$ is purely transcendental, so we may consider $L'^{im \beta_i}$ as subfield of $F_i^\star$ and denote the induced valuation on $F'$ by $v'_i$. Now let $F'_i$ be the relative algebraic closure of $F'$ in $F_i^\star$ and observe that then $F_i \subseteq F'_i$, that $F_i$ is existentially closed in $F'_i$ (because it is in $F_i^\star$) and that $L' \cap F_i^\star = L'^{im \beta_i}$. Therefore, $res: G_{F'} \to A = Gal(L'/F')$ gives rise to a solution of our lifted locally split embedding problem, which, in particular, is a locally conjugate solution. □

# References

[En] O.Endler: *Valuation theory*, Springer (1972).

[Er1] Y.Ershov: *Semilocal fields*, Soviet Math.Dokl. **15** (1974), 424-428.

[Er2] Y.Ershov: *Projectivity of absolute Galois groups of $RC_\zeta^\star$-fields*, Proc. III Int. Conf.Algebra Krasnoyarsk 93, deGruyter (1996), 63-80.

[Er3] Y.Ershov: *Free products of absolute Galois groups*, Dokl.Math. **56(3)** (1997), 915-917.

[Ge] W.-D.Geyer: *Galois groups of intersections of local fields*, Isr.J.Math. **30** (1978), 382-396.

[Gr] K.W. Gruenberg: *Projective profinite groups*, Journal London Math.Soc. **42** (1967) 155-165.

[Ha] D.Haran: *On closed subgroups of free products of profinite groups*, Proc.London Math.Soc. **55** (1987), 266-298.

[HaJ] D.Haran, M.Jarden: *The absolute Galois group of a pseudo p-adically closed field*, J.Reine Angew.Math. **383** (1988), 147-206.

[Hs] H.Hasse: *Existenz und Mannigfaltigkeit abelscher Algebren mit vorgegebener Galoisgruppe über einem Teilkörper des Grundkörpers I*, Math.Nachrichten **1** (1948), 40-61.

[He] B.Heinemann: *On finite intersections of 'henselian valued' fields*, manuscr.math. **52** (1985), 37-61.

[HP] B.Heinemann, A.Prestel: *Fields regularly closed with respect to finitely many valuations and orderings*, Canad.Math.Soc.Conf.Proc. **4** (1984), 297-336.

[J] M.Jarden: *Infinite Galois theory*, in: M.Hazewinkel (ed.): Handbook of Algebra, vol. 1 (1996).

[K] F.-V.Kuhlmann: *Valuation theory of fields, abelian groups and modules*, to appear.

[M] O.V.Mel'nikov: *On free products of absolute Galois groups*, (1997, Russian), English transl. in Siberian Math.J. **40 (1)** (1999), 95-99.

[Po1] F.Pop: *Classically projective groups and pseudo classically closed fields*, Prerpint Math.Inst.Heidelberg 1990.

[Po2] F.Pop: *Embedding problems over large fields*, Annals of Math. **144** (1996), 1-34.

[PrC] S.Prieß-Crampe *Angeordnete Strukturen: Gruppen, Körper, projektive Ebenen*, Springer Erg.d.Math **98** (1983).

[Ri] P.Ribenboim: *Théorie des valuations*, Montreal (1968).

[Ru] R.Rumely: *Arithmetic over the ring of all algebraic integers*, J.reine angew.Math. **368** (1986), 127-133.

*e-mail:* `Jochen.Koenigsmann@uni-konstanz.de`